\documentclass[11pt,a4paper]{amsart}
\usepackage[T1]{fontenc}
\usepackage{lmodern}
\usepackage[margin=2.7cm]{geometry}
\usepackage{amsmath,amssymb,amsthm,mathtools}
\usepackage{enumitem}
\usepackage{microtype}
\usepackage{xcolor}
\usepackage[colorlinks,linkcolor=blue!50!black,citecolor=blue!50!black,urlcolor=blue!50!black]{hyperref}
\hypersetup{pdftitle={Randomized versus deterministic optimization for convex and symmetric input classes},
 pdfauthor={David Krieg}}

\setlist{leftmargin=2.2em}
\numberwithin{equation}{section}

\newtheorem{theorem}{Theorem}[section]
\newtheorem{proposition}[theorem]{Proposition}
\newtheorem{lemma}[theorem]{Lemma}
\newtheorem{corollary}[theorem]{Corollary}
\newtheorem{conjecture}[theorem]{Conjecture}
\theoremstyle{definition}
\newtheorem{problem}[theorem]{Open Problem}
\theoremstyle{remark}
\newtheorem{remark}[theorem]{Remark}

\newcommand{\R}{\mathbb{R}}
\newcommand{\N}{\mathbb{N}}
\newcommand{\E}{\mathbb{E}}
\renewcommand{\P}{\mathbb{P}}
\newcommand{\dt}{\mathrm{det}}
\newcommand{\ran}{\mathrm{ran}}
\newcommand{\eps}{\varepsilon}
\newcommand{\vol}{\operatorname{vol}}
\newcommand{\conv}{\operatorname{conv}}
\newcommand{\spn}{\operatorname{span}}

\newcommand{\dist}{\operatorname{dist}}
\newcommand{\Nf}{\mathcal N}
\newcommand{\OPT}{\mathrm{OPT}}
\newcommand{\APP}{\mathrm{APP}_\infty}
\newcommand{\abs}[1]{\lvert#1\rvert}
\newcommand{\norm}[1]{\lVert#1\rVert}
\newcommand{\ip}[2]{\langle#1,#2\rangle}

\title[Randomized versus deterministic optimization]{Randomized versus deterministic optimization\\ for convex and symmetric input classes}
\author{David Krieg}
\address{University of Passau, Germany}
\email{david.krieg@uni-passau.de}
\date{}

\begin{document}

\begin{abstract}
We study global optimization with function evaluations only: we want to compute the supremum of a function from a
convex and symmetric class $F$ of bounded functions, using $n$ function values, that is, $n$ queries to a
zeroth-order oracle. We compare deterministic algorithms with randomized
algorithms in the worst case over $F$. We show that for every such class the error of the best deterministic algorithm with
$8n-1$ queries is at most $Cn^2$ times the error of the best randomized algorithm with
$n$ queries, where $C$ is a universal constant. The same holds, up to one function value, for the problems of finding a maximizer or minimizer. 
On the other hand, for all $m$ and $n$ there is a class for which deterministic algorithms
with $m$ function values are worse than randomized algorithms with $n$ function values by a factor of
order $n^{3/2}$, up to a logarithmic factor, for the computation of the supremum. The upper bound is obtained from a lower bound for randomized algorithms in
terms of Bernstein widths. We conjecture that this bound can be improved by a factor $\sqrt{n}$, which would make the
exponent $3/2$ sharp. We also show that, for randomized algorithms, finding a point at which the supremum is almost
attained can be harder than computing the supremum by a factor of order $n/\log n$, while the two problems are
known to be equivalent for deterministic algorithms and convex and symmetric input classes.
\end{abstract}


\maketitle

\section{Introduction and results}

Let $D$ be a nonempty set and let $B(D)$ be the space of bounded real-valued functions on $D$ with the norm
$\norm{f}_\infty=\sup_{x\in D}\abs{f(x)}$. We want to compute the supremum
\[
 S(f)\,:=\,\sup_{x\in D} f(x)
\]
of an unknown function $f$ from a class $F\subset B(D)$ with the help of $n$ function values. Throughout,
$F$ is assumed to be nonempty, convex and symmetric, that is, $-f\in F$ whenever $f\in F$. Here, convexity is a
property of the class and not of the functions in it, which need not be convex or concave: typical examples are the
unit balls of Lipschitz, H\"older or Sobolev spaces and of reproducing kernel Hilbert spaces, provided that these
consist of bounded functions. Since $F$ is symmetric, it makes no difference whether we maximize or minimize. This is the
problem of global optimization in the worst-case setting of information-based complexity; we refer to
\cite{NY,No88,TWW} for background. In the language of optimization, an algorithm has access to $f$ only through a
zeroth-order (or black-box) oracle, which returns the exact value $f(x)$ for a query $x\in D$, and $n$ is the
number of queries. The minimal number of queries for a given error is the (deterministic or
stochastic) complexity of the class in the sense of Nemirovsky and Yudin \cite[Section~1.4]{NY}. The functions in $F$ are bounded, but we do not assume that $F$ is a bounded
subset of $B(D)$. Therefore the error numbers and widths that are defined below may be infinite, and all
inequalities between them are understood in $[0,\infty]$; statements about ratios are to be read as the corresponding inequalities. Throughout, $\log$ denotes the natural logarithm.

An algorithm may choose the points at which $f$ is evaluated one after the other, depending on the function values
that have already been observed, that is, adaptively or sequentially. A deterministic algorithm does this according to a fixed rule, while a randomized
algorithm may also use random numbers.
We write $e_n^{\dt}(S,F)$ for the minimal worst-case error, or minimax error, of deterministic
algorithms that use $n$ function values, and $e_n^{\ran}(S,F)$ for the
minimal error of randomized algorithms, where the error of a randomized algorithm is its expected error for
the worst function in $F$. The precise definitions are given in Section~\ref{sec:prelim}.
Besides the optimal value $S(f)$, one may ask for a point $\hat x\in D$ at which it is almost attained. We denote
this problem by $\OPT$. Here the error is the residual $S(f)-f(\hat x)$, which is also called the optimality gap or
the simple regret, and the minimal errors $e_n^{\dt}(\OPT,F)$ and $e_n^{\ran}(\OPT,F)$ are defined accordingly, see
Section~\ref{sec:opt}.
The question of this paper is:
\begin{quote}
\emph{How much better than deterministic algorithms can randomized algorithms be?}
\end{quote}
They are never worse, up to a factor of two, see Proposition~\ref{prop:meas},
but they might be significantly better.
For deterministic algorithms it is known that the squential choice of the points does not help much: on convex and
symmetric classes one loses at most a factor two if all points are chosen in advance, that is, nonadaptively or in one
batch, see
Proposition~\ref{prop:det} below, which goes back to
Wasilkowski~\cite{Wa84} and Novak~\cite{No88,No96}; for Lipschitz classes see already Sukharev~\cite{Su71}. So the
question is about randomization, and all our upper bounds
for deterministic algorithms hold for algorithms that choose their points in advance. It is also known that
randomization does not help for many
classical classes, like H\"older classes \cite{NY,No88}. For general convex and symmetric classes, however, no bound
on the gain of randomization seems to be known.
Our first result is that the gain is at most polynomial in $n$.

\begin{theorem}\label{thm:upper}
Let $F\subset B(D)$ be nonempty, convex and symmetric. Then, for all $n\ge 50$,
\[
 e_{8n-1}^{\dt}(S,F)\;\le\; C\,n^{2}\; e_n^{\ran}(S,F)
\]
with the absolute constant $C=22\,400$.
\end{theorem}

The same bound holds for the problem of finding a maximizer, up to a factor two and one function value:
by Corollary~\ref{cor:opt-upper}, for all
$n\ge50$,
\[
 e_{8n-1}^{\dt}(\OPT,F)\;\le\;2\,C\,n^{2}\;e_{n-1}^{\ran}(\OPT,F).
\]

Our second result shows that the gain of randomization can indeed be large, and larger than for linear problems with arbitrary continuous linear measurements, where it is at
most of order $n$, see \cite{KNU,KNU2,Ul26} and the examples in \cite{KNW,KW}. Moreover, the large loss of deterministic algorithms cannot be
compensated by giving the deterministic algorithm more function values.

\begin{theorem}\label{thm:lower}
There is an absolute constant $c>0$ such that for all $m,n\in\N$ there is a set $D$ and a nonempty, bounded,
convex and symmetric class $F\subset B(D)$ with
\[
 e_{m}^{\dt}(S,F)\;\ge\; c\,\frac{n^{3/2}}{\log(n+1)}\; e_n^{\ran}(S,F).
\]
More precisely, $e_m^{\dt}(S,F)\ge\frac18$ and $e_n^{\ran}(S,F)\le C\,n^{-3/2}\log(n+1)$ for an absolute
constant~$C$.
\end{theorem}

The class $F$ in
Theorem~\ref{thm:lower} comes from a problem that may be of independent interest, the estimation of
the norm $\max\{\norm{\xi}_\infty,\tau\norm{\xi}_1\}$ of a vector $\xi\in\R^d$ with $\norm\xi_1\le1$ from $n$ linear
measurements, where the dimension $d$ is larger than $6^m$ and $\tau$ is of order $\log(n)/n$. Every problem of this
kind, with an arbitrary norm on a convex and symmetric set, is the problem of computing the supremum for a class
of linear functions on the unit ball of the dual norm, see Lemma~\ref{lem:dict}. A deterministic
algorithm cannot find a large coordinate of $\xi$ among exponentially many, while a randomized algorithm
can find all coordinates that matter and estimates the remaining term $\tau\norm{\xi}_1$ by a Monte Carlo method.

Theorem~\ref{thm:lower} is a statement about the problem $S$ of finding the optimal value. For the problem of finding a maximizer, it follows from
a result of Kunsch and Wnuk~\cite{KW} that a gain of order $n/\log n$ is possible if $m$ is at most polynomial in $n$, see
Proposition~\ref{prop:opt-lower}; we
do not know whether a gain of order $n^{3/2}$ is possible, see Open Problem~\ref{prob:opt-gain}.
For deterministic algorithms, computing the supremum and finding a maximizer are equivalent, up to a factor of four
in the error, see Propositions~\ref{prop:det} and~\ref{prop:opt}. For randomized algorithms this is not the case
and it is an interesting problem on its own to determine
how much harder it can be to find a maximizer than to find the maximum 
if randomized algorithms are allowed.
Theorem~\ref{thm:point-value} shows that
there are classes for which finding a maximizer with $n$ function values is harder than computing the supremum by a
factor of order $n/\log n$,
while the ratio is never larger than $Cn^2$, 
at least if the number of function values for the harder problem $\OPT$ is increased by a constant factor.

So the largest possible gain of randomization for the computation of the supremum, measured by the ratio $e^{\dt}_m/e^{\ran}_n$ with
$m\ge 8n$, is between $n^{3/2}$ and $n^2$, up to a logarithmic factor. We believe that $n^{3/2}$ is the truth. To explain
why, we describe how Theorem~\ref{thm:upper} is proved. For $k\in\N_0$, the $k$-th \emph{Bernstein width} of $F$
in $B(D)$ is
\[
 b_k(F)\,:=\,\sup\bigl\{\lambda\ge0\colon\ \lambda\,(B_{B(D)}\cap V)\subset F\ \text{ for some subspace }
 V\subset B(D),\ \dim V=k+1\bigr\},
\]
where $B_{B(D)}$ is the unit ball of $B(D)$ and $b_k(F):=0$ if there is no such subspace. In other words, $b_k(F)$ is
the largest radius of a $(k+1)$-dimensional ball in the uniform norm that is contained in $F$.
Bernstein widths are a classical tool for lower bounds: if $F$ contains a ball of dimension larger than $n$, then
there is a nonzero function in this ball that vanishes at any given $n$ points, and a deterministic algorithm cannot
distinguish it from the zero function. We show that Bernstein widths bound the error of deterministic and
randomized algorithms from above and from below, respectively.

\begin{theorem}\label{thm:bernstein}
Let $F\subset B(D)$ be nonempty, convex and symmetric. Then
\begin{align}
 e_{8n-1}^{\dt}(S,F)\;&\le\; 32\,n\; b_{4n-1}(F) &&\text{for all } n\in\N, \label{eq:det-bern}\\
 e_{n}^{\ran}(S,F)\;&\ge\; \frac{c_0}{4n}\; b_{4n-1}(F) &&\text{for all } n\ge 50, \label{eq:ran-bern}
\end{align}
where $c_0=1/175$. The lower bound \eqref{eq:ran-bern} also holds for randomized algorithms that use $n$ arbitrary
linear measurements of $f$ instead of $n$ function values, that is, for the number
$e_n^{\ran}(S,F,\Lambda^{\rm all})$ defined in Section~\ref{sec:prelim}.
\end{theorem}

Theorem~\ref{thm:upper} is an immediate consequence, and it also holds with $e_n^{\ran}(S,F,\Lambda^{\rm all})$ on
the right-hand side. In particular, the bound of Theorem~\ref{thm:upper} remains valid if both algorithms may use further linear
measurements besides function values, for example values of directional derivatives, as provided by a first-order
oracle. Here, every real number that is observed counts as one measurement, so that a gradient in dimension $d$
counts as $d$ measurements. The lower bound of Theorem~\ref{thm:lower} also remains valid in this setting, since
the deterministic bound for the classes considered there holds for arbitrary linear measurements, see
Lemma~\ref{lem:dict} and Proposition~\ref{prop:ex-det}. 
The factor $n$ in \eqref{eq:det-bern} cannot be improved
in general, see Remark~\ref{rem:sharp}; in fact, for the classes of Theorem~\ref{thm:lower} with $m\ge8n$ both sides
of \eqref{eq:det-bern} are of the same order up to a logarithm, see Remark~\ref{rem:example}(2). The factor $1/n$ in
\eqref{eq:ran-bern}, on the other hand, is probably not optimal. Let us note that for linear problems Bernstein
widths are lower bounds for the error of randomized algorithms without any loss in $n$, up to
constants, see \cite{He92,No92,Ku16,Ku17,KNU}. The proof of this fact uses the linearity of the problem, and the
statement is not true for the functional $S$, as the example after the following conjecture shows.

\begin{conjecture}\label{conj:main}
There are absolute constants $c>0$ and $C\in\N$ such that for every nonempty, convex and symmetric class
$F\subset B(D)$ and all $n\in\N$,
\[
 e_{n}^{\ran}(S,F)\;\ge\; \frac{c}{\sqrt n}\; b_{Cn}(F).
\]
\end{conjecture}

Together with the deterministic bound (in the form of Proposition~\ref{prop:rho-bern}), the conjecture would imply
$e^{\dt}_{2Cn-1}(S,F)\le C' n^{3/2}\,e^{\ran}_n(S,F)$ for all classes, and by Theorem~\ref{thm:lower} this
would be optimal up to a logarithmic factor. The factor $n^{-1/2}$ in the conjecture cannot be improved: if $F$ is
the class of all linear functions $x\mapsto\ip x\xi$ with $\norm\xi_1\le1$ on the cube $[-1,1]^d$, then
$S(f)=\norm{f}_\infty=\norm{\xi}_1$, so that $F$ is itself a ball and $b_{k}(F)=1$ for all $k<d$, while
$e_n^{\ran}(S,F)\le 10\,n^{-1/2}$ for all $n$ and $d$, see Remark~\ref{rem:norms}(2). For the purpose of comparing
$e^{\dt}$ and $e^{\ran}$ a version of the conjecture with additional logarithmic factors, or with $b_{m(n)}$ for
some $m(n)=n^{1+o(1)}$, would be just as interesting.

\subsection*{Remarks}
(1) \emph{The classes in Theorem~\ref{thm:lower} depend on $m$ and $n$.} Theorem~\ref{thm:lower} is a statement
about the errors for given numbers of function values, like Theorem~\ref{thm:upper}, and not about rates of
convergence for a fixed class. For a fixed class the maximal difference of the polynomial rates of $e_n^{\dt}$ and
$e_n^{\ran}$ is a different question, see Open Problem~\ref{prob:rates}.

(2) \emph{Oversampling.} In Theorem~\ref{thm:upper} the deterministic algorithm uses about eight times as many function
values as the randomized one. With four times as many, there is still a polynomial bound, namely
$e^{\dt}_{4n-1}(S,F)\le C\,n^3\,e^{\ran}_n(S,F)$ for $n\ge50$, with $C$ as in Theorem~\ref{thm:upper}; this follows from Proposition~\ref{prop:rho-bern} with $k=1$
and \eqref{eq:ran-bern}. We do not know whether a
polynomial bound also holds without oversampling, that is, for $e^{\dt}_n/e^{\ran}_n$, see Open Problem~\ref{prob:oversampling}.

(3) \emph{Other problems.} All upper bounds for deterministic algorithms in this paper also hold for the problem of
uniform approximation and for the computation of the norm $\norm{f}_\infty$ instead of $S(f)$, see
Proposition~\ref{prop:det} and the remark after it; the lower bounds for randomized algorithms are proved for the norm first and then
transferred to $S$.

(4) \emph{Unbounded classes.} The class $F$ need not be bounded. An example is the class of all functions on
$D=[0,1]$ with Lipschitz constant at most one, which contains all constant functions; which contains all constant functions; here $e_n^{\dt}(S,F)=1/(4n)$. For unbounded classes the
errors may be infinite, and Theorem~\ref{thm:upper} shows in particular that $e^{\dt}_{8n-1}(S,F)$ is finite whenever
$e^{\ran}_n(S,F)$ is finite and $n\ge50$. The classes in Theorem~\ref{thm:lower} are bounded.

(5) \emph{Symmetry.} The symmetry of $F$ is essential. For classes $F\subset B(D)$ that are convex but not symmetric, 
already for deterministic algorithms,
the adaptive choice of the points can help a lot, see \cite{No96,KNU2},
and optimization can be a lot easier than uniform approximation.
Moreover, the gain of randomization can be
unbounded: an extreme example is the convex class $F$ of all functions $f\colon[0,1]\to[0,1]$ for which there is a
constant $c$ with $f\le c$ and $f(x)=c$ for all but finitely many $x$; here $e_n^{\dt}(S,F)=\frac12$ for all $n$, while 
$e_1^{\ran}(S,F)=0$. 

\subsection*{Organization} Section~\ref{sec:prelim} contains the definitions, some elementary reductions, the relation to the estimation of
norms from linear measurements (Lemma~\ref{lem:dict}) and a remark on linear algorithms and measurability.
Sections~\ref{sec:det} and~\ref{sec:ran} prove the two parts of Theorem~\ref{thm:bernstein}, and hence
Theorem~\ref{thm:upper}. Section~\ref{sec:example} proves Theorem~\ref{thm:lower}; it can be read independently of
the proofs in Sections~\ref{sec:det} and~\ref{sec:ran}. Section~\ref{sec:opt} contains the consequences
for the problem of finding a maximizer, and open problems are collected in
Section~\ref{sec:open}.

\section{Preliminaries}\label{sec:prelim}

\subsection{Algorithms and errors}
Besides $S$ we consider the functional $\Nf(f):=\norm{f}_\infty$. Note that $\Nf(f)=\max\{S(f),S(-f)\}$.

A \emph{deterministic algorithm} using $n$ function values is a mapping $A_n\colon F\to\R$ of the form
\[
 A_n(f)=\varphi(y_1,\dots,y_n),\qquad y_i=f(x_i),\quad x_i=\psi_i(y_1,\dots,y_{i-1}),
\]
with arbitrary mappings $\psi_i\colon\R^{i-1}\to D$ and $\varphi\colon\R^n\to\R$. The algorithm is called
\emph{nonadaptive} if all $\psi_i$ are constant, that is, if the points $x_1,\dots,x_n$ do not depend on the previously observed values.
For $T\in\{S,\Nf\}$ the worst-case error of $A_n$ is $e(A_n,T,F):=\sup_{f\in F}\abs{T(f)-A_n(f)}$, and we put
\[
 e_n^{\dt}(T,F):=\inf_{A_n}\,e(A_n,T,F),
\]
where the infimum is taken over all deterministic algorithms using $n$ function values.
A \emph{randomized algorithm} is a family $(A_n^\omega)_{\omega\in\Omega}$ of deterministic algorithms, indexed by a
probability space $(\Omega,\mathcal A,\P)$. It is called nonadaptive if all $A_n^\omega$ are nonadaptive. Its error is
\[
 e\bigl((A_n^\omega),T,F\bigr):=\sup_{f\in F}\,\E_\omega\abs{T(f)-A_n^\omega(f)} ,
\]
and $e_n^{\ran}(T,F)$ is the infimum of this error over all measurable randomized algorithms. Similarly to
\cite{KNU2}, measurability means that for every finite-dimensional subspace of $B(D)$ that is spanned by functions from $F$, with basis
$g_1,\dots,g_r$, the mappings $(\omega,y)\mapsto\bigl(g_l(\psi_i^\omega(y))\bigr)_{l\le r}$ and
$(\omega,y)\mapsto\varphi^\omega(y)$ are measurable with respect to the product of $\mathcal A$ and the Borel
$\sigma$-algebra. Then, for inputs $f=\sum_l\xi_lg_l$ from this subspace, the
observed function values and the output are jointly measurable in $\omega$ and $\xi$. This is what is needed for the averaging argument of Bakhvalov~\cite{Ba59}, which we use in the following
form: if $\mu$ is a Borel measure on such a subspace that is concentrated on a Borel subset of $F$ and has total
mass at most one, then
\begin{equation}\label{eq:bakhvalov}
 \sup_{f\in F}\E_\omega\abs{T(f)-A_n^\omega(f)}\;\ge\;\E_\omega\int\abs{T(f)-A_n^\omega(f)}\,d\mu(f)
 \;\ge\;\inf_{\omega}\int\abs{T(f)-A_n^\omega(f)}\,d\mu(f).
\end{equation}
So lower bounds for randomized algorithms follow from lower bounds for the average error of deterministic algorithms
with respect to $\mu$. In the following, all randomized algorithms are assumed to be measurable.

Our lower bounds for randomized algorithms do not use that the measurements are function values. Let
$\Lambda^{\rm all}$ be the set of all linear functionals on the linear span of $F$; they need not be continuous. An
algorithm using $n$ linear measurements is defined as above, with $y_i=L_i(f)$ and
$L_i=\psi_i(y_1,\dots,y_{i-1})\in\Lambda^{\rm all}$ in place of $y_i=f(x_i)$. For randomized algorithms, measurability
now means that the mappings $(\omega,y)\mapsto\bigl(\psi_i^\omega(y)(g_l)\bigr)_{l\le r}$ and
$(\omega,y)\mapsto\varphi^\omega(y)$ are measurable for the subspaces and bases above. Since
$L(\sum_l\xi_lg_l)=\sum_l\xi_lL(g_l)$, the observed values and the output are again jointly measurable in $\omega$ and
$\xi$, by induction on the number of measurements, and \eqref{eq:bakhvalov} holds as before. We write $e_n^{\ran}(T,F,\Lambda^{\rm all})$ for the infimum of the error over all measurable randomized
algorithms using $n$ linear measurements. Function values are linear measurements, and hence
$e_n^{\ran}(T,F,\Lambda^{\rm all})\le e_n^{\ran}(T,F)$. When we want to emphasize that algorithms use function
values, we add $\Lambda^{\rm std}$ to the notation, for example $e_n^{\dt}(T,F,\Lambda^{\rm std}):=e_n^{\dt}(T,F)$.

\begin{remark}\label{rem:ran-det}
We do not impose any measurability on deterministic algorithms. Therefore a deterministic algorithm is not
necessarily a randomized algorithm in the above sense, and the inequality $e_n^{\ran}(T,F)\le e_n^{\dt}(T,F)$, which
one would expect, is not clear in general. For convex and symmetric classes it holds up to a constant, see
Proposition~\ref{prop:meas} and Proposition~\ref{prop:opt}.
\end{remark}

\subsection{The radius of information}
For points $x_1,\dots,x_n\in D$ let
\begin{align*}
 \rho(x_1,\dots,x_n;F)&:=\sup\bigl\{\norm{h}_\infty\colon h\in F,\ h(x_1)=\dots=h(x_n)=0\bigr\},\\
 \rho_n(F)&:=\inf_{x_1,\dots,x_n\in D}\rho(x_1,\dots,x_n;F).
\end{align*}
The number $\rho_n(F)$ is the smallest possible radius of the set of functions in $F$ that cannot be distinguished
from zero by $n$ function values (it may be infinite if $F$ is unbounded). It is nonincreasing in $n$, and $\rho_n(F)\ge b_n(F)$ since every
$(n+1)$-dimensional space contains a nonzero function that vanishes at $n$ given points.

The following proposition is due to Wasilkowski~\cite[Theorems~4.1 and~5.1]{Wa84}, who considers the norm on convex
and symmetric subsets of a normed space and the minimum of continuous functions on $[0,1]$, with measurements
given by arbitrary linear functionals. In the present setting of function values and bounded functions on an
arbitrary set, the statement for $S$ is contained in \cite[Propositions~1.2.6 and~1.3.3]{No88}. We include the short
proof.

\begin{proposition}\label{prop:det}
Let $F\subset B(D)$ be nonempty, convex and symmetric, and $n\in\N$. Then
\[
 \tfrac12\,\rho_n(F)\;\le\; e_n^{\dt}(S,F)\;\le\;\rho_n(F).
\]
The upper bound $e_n^{\dt}(T,F)\le\rho_n(F)$ also holds for $T=\Nf$, and for $T\in\{S,\Nf\}$ it is attained, up to
any $\eta>0$, by nonadaptive algorithms.
\end{proposition}

\begin{proof}
\emph{Lower bound.} Let $A_n$ be a deterministic algorithm and let $x_1,\dots,x_n$ be the points that it uses
for the input $f=0$, which belongs to $F$. Put $y_0:=A_n(0)$. If $h\in F$ vanishes at $x_1,\dots,x_n$, then the
algorithm observes the same values for $h$, for $-h$ and for $0$, hence it uses the same points and returns
$A_n(h)=A_n(-h)=y_0$. Since $S(0)=0$, the error of $A_n$ is at least
\[
 \max\bigl\{\abs{y_0},\,\abs{S(h)-y_0}\bigr\}\ge\tfrac12 S(h)
 \qquad\text{and}\qquad
 \max\bigl\{\abs{y_0},\,\abs{S(-h)-y_0}\bigr\}\ge\tfrac12 S(-h).
\]
As $\max\{S(h),S(-h)\}=\norm{h}_\infty$, we get $e(A_n,S,F)\ge\frac12\norm h_\infty$. Taking the supremum over $h$
gives $e(A_n,S,F)\ge\frac12\rho(x_1,\dots,x_n;F)\ge\frac12\rho_n(F)$.

\smallskip
\emph{Upper bound.} We may assume that $\rho_n(F)$ is finite. Fix $x_1,\dots,x_n\in D$ such that
$\rho:=\rho(x_1,\dots,x_n;F)$ is finite. For $y\in\R^n$ let
$F_y:=\{g\in F\colon g(x_i)=y_i \text{ for } i\le n\}$. If $f,g\in F_y$, then $\frac12(f-g)\in F$ by convexity and
symmetry, and this function vanishes at the points, so that $\norm{f-g}_\infty\le2\rho$. Since
$\abs{T(f)-T(g)}\le\norm{f-g}_\infty$ for $T\in\{S,\Nf\}$, the set $\{T(g)\colon g\in F_y\}$ is contained in an
interval of length $2\rho$. The nonadaptive algorithm that returns the midpoint of this interval for
$y=(f(x_1),\dots,f(x_n))$ has error at most $\rho$. Taking the infimum over all such points proves the claim.
\end{proof}

\begin{remark}\label{rem:app}
The upper bound also holds for the problem of uniform approximation. Let $x_1,\dots,x_n\in D$ be such that
$\rho(x_1,\dots,x_n;F)$ is finite. Then the nonadaptive algorithm that returns, for
$y=(f(x_1),\dots,f(x_n))$, the function $x\mapsto\frac12\bigl(\sup_{g\in F_y}g(x)+\inf_{g\in F_y}g(x)\bigr)$ is well
defined, since $\norm{f-g}_\infty\le2\rho(x_1,\dots,x_n;F)$ for $f,g\in F_y$, and it
approximates every $f\in F$ with an error of at most $\rho(x_1,\dots,x_n;F)$ in the uniform norm. Hence all upper
bounds for deterministic algorithms in this paper, which are obtained from bounds on $\rho_n(F)$, hold for uniform
approximation as well. We write $e_n^{\dt}(\APP,F)$ for the infimum of $\sup_{f\in F}\norm{f-A_n(f)}_\infty$ over all
deterministic algorithms $A_n\colon F\to B(D)$, which are defined as above with $\varphi\colon\R^n\to B(D)$. Then
\[
 e_n^{\dt}(\APP,F)=\rho_n(F),
\]
see \cite[Proposition~1.2.6]{No88}. The upper bound was just shown. For the lower bound let $x_1,\dots,x_n$ be the points used for the input $f=0$ and
let $h\in F$ vanish at these points. As in the proof of Proposition~\ref{prop:det}, the algorithm returns the same
function $g$ for $h$ and for $-h$, and $\norm{h-g}_\infty+\norm{-h-g}_\infty\ge2\norm h_\infty$.

Moreover, linear algorithms are optimal for uniform approximation on convex and symmetric classes: the number
$e_n^{\dt}(\APP,F)$ does not change if the infimum is restricted to algorithms of the form
$A_n(f)=\sum_{i\le n}f(x_i)\,a_i$ with $x_i\in D$ and $a_i\in B(D)$, see \cite[Proposition~1.2.6]{No88} and Lemma~\ref{lem:linear}. So all upper
bounds for $\rho_n(F)$ in this paper are also bounds for the minimal error of linear algorithms.
\end{remark}

\begin{remark}[Adaption]\label{rem:adaption}
The upper bounds in Proposition~\ref{prop:det} and Remark~\ref{rem:app} are obtained with nonadaptive algorithms. So,
for deterministic algorithms on convex and symmetric classes, the adaptive choice of the points helps at most by a
factor of two for $S$, and it does not help at all for uniform approximation. This is the reason why we do not
distinguish between adaptive and nonadaptive deterministic algorithms in this paper: all our upper bounds for
$e_n^{\dt}$ are bounds for $\rho_n(F)$ and therefore hold for nonadaptive algorithms.
\end{remark}

\subsection{From the supremum to the norm}
Lower bounds for randomized algorithms are more convenient for the norm $\Nf$. The following lemma transfers them
to $S$ at the price of a factor two in the number of function values.

\begin{lemma}\label{lem:StoN}
For every nonempty symmetric $F\subset B(D)$ and $n\in\N$,
\[
 e_{2n}^{\ran}(\Nf,F)\le2\,e_n^{\ran}(S,F)\qquad\text{and}\qquad
 e_{2n}^{\ran}(\Nf,F,\Lambda^{\rm all})\le2\,e_n^{\ran}(S,F,\Lambda^{\rm all}).
\]
\end{lemma}

\begin{proof}
Let $(A^\omega)$ be a measurable randomized algorithm for $S$ using $n$ function values, and define
$\tilde A^{\omega}(f):=\max\{A^\omega(f),A^{\omega}(-f)\}$. A function value of $-f$
is obtained from a function value of $f$, so $\tilde A^\omega$ uses $2n$ function values. The same holds for
linear measurements, since $L(-f)=-L(f)$ for $L\in\Lambda^{\rm all}$. In both cases the decision rules and the output of $\tilde A^\omega$ are obtained from
those of $A^\omega$ by a change of sign of the arguments and a maximum, so that $(\tilde A^\omega)$ is measurable.
Since $\norm f_\infty=\max\{S(f),S(-f)\}$ and $\abs{\max\{a,b\}-\max\{a',b'\}}\le\abs{a-a'}+\abs{b-b'}$,
\[
 \E_\omega\abs{\Nf(f)-\tilde A^\omega(f)}\le\E_\omega\abs{S(f)-A^\omega(f)}+\E_\omega\abs{S(-f)-A^{\omega}(-f)}\le2\,e\bigl((A^\omega),S,F\bigr)
\]
for every $f\in F$, because $-f\in F$.
\end{proof}

\subsection{Classes of linear functions and the estimation of norms}
We also consider the problem of estimating an arbitrary norm on an arbitrary convex and symmetric subset of $\R^d$
from $n$ linear measurements. The following lemma shows that this is, in fact,
merely a special case of the optimization problems considered in this paper. This will be used in
Section~\ref{sec:example} to construct an example where randomized algorithms are much better than deterministic
algorithms. Let $p$ be a norm on $\R^d$ and let $B\subset\R^d$. Algorithms for the approximation of $p$ on $B$
using $n$ linear measurements are defined as above, with the input $\xi\in B$, with
$y_i=\ip{a_i}{\xi}$ and $a_i=\psi_i(y_1,\dots,y_{i-1})\in\R^d$ in place of $y_i=f(x_i)$, and with the error
$\abs{p(\xi)-A_n(\xi)}$. A randomized algorithm of this kind is called measurable if the mappings
$(\omega,y)\mapsto\psi_i^\omega(y)$ and $(\omega,y)\mapsto\varphi^\omega(y)$ are measurable. We write
$e_n^{\dt}(p,B,\Lambda^{\rm all})$ and $e_n^{\ran}(p,B,\Lambda^{\rm all})$ for the corresponding minimal errors; as before,
$\Lambda^{\rm all}$ indicates that all linear measurements are admitted.

\begin{lemma}\label{lem:dict}
Let $p$ be a norm on $\R^d$ with unit ball $B_p$, let $B\subset\R^d$ be convex and symmetric with $\spn B=\R^d$, and let
\[
 K:=\bigl\{x\in\R^d\colon\ip x\xi\le p(\xi)\text{ for all }\xi\in\R^d\bigr\}
 \qquad\text{and}\qquad
 F:=\bigl\{f_\xi\colon\xi\in B\bigr\}\subset B(K),
\]
where $f_\xi(x):=\ip x\xi$ for $x\in K$. Then $F$ is nonempty, convex and symmetric,
$S(f_\xi)=\norm{f_\xi}_\infty=p(\xi)$ for all $\xi\in\R^d$, and for all $n\in\N$,
\[
 e_n^{\dt}(S,F)=e_n^{\dt}(p,B,\Lambda^{\rm all})\qquad\text{and}\qquad e_n^{\ran}(S,F)=e_n^{\ran}(S,F,\Lambda^{\rm all})=e_n^{\ran}(p,B,\Lambda^{\rm all}).
\]
The same holds if only nonadaptive algorithms are admitted. Moreover, for all $k\in\N_0$,
\[
 b_k(F)=\sup\bigl\{\lambda\ge0\colon\lambda\,(B_p\cap U)\subset B\ \text{ for some subspace }U\subset\R^d,\ \dim U=k+1\bigr\},
\]
where the supremum of the empty set is zero. In particular, $b_k(F)=1$ for $k<d$ if $B=B_p$.
\end{lemma}

\begin{proof}
The set $K$ is the unit ball of the dual norm $p^*(a):=\max\{\ip a\xi\colon p(\xi)\le1\}$. It is a symmetric convex
body, and $\sup_{x\in K}\ip x\xi=p(\xi)$ for all $\xi$ by the Hahn--Banach theorem. Since $K=-K$, this gives
$S(f_\xi)=\norm{f_\xi}_\infty=p(\xi)$, and the mapping $\xi\mapsto f_\xi$ is linear and injective, so that $F$ is
nonempty, convex and symmetric.

\smallskip
\emph{Algorithms.} The function value of $f_\xi$ at $x\in K$ is the linear measurement $\ip x\xi$. On the other
hand, for $a\in\R^d\setminus\{0\}$ the point $x_a:=a/p^*(a)$ belongs to $K$, and $\ip a\xi=p^*(a)\,f_\xi(x_a)$; for
$a=0$ this also holds if $x_0$ is any fixed point of $K$. Hence every algorithm for $S$ on $F$ using $n$ function
values is an algorithm for $p$ on $B$ using $n$ linear measurements, and every algorithm of the second kind is
turned into one of the first kind, with the same output for every input, by evaluating $f_\xi$ at $x_{a_i}$
instead of measuring $\ip{a_i}\xi$ and multiplying the observed value by $p^*(a_i)$. Nonadaptive algorithms stay
nonadaptive. Since $\spn B=\R^d$, the linear span of $F$ consists of all functions $f_\xi$, $\xi\in\R^d$, and with
the basis $g_l:=f_{e_l}$, $l\le d$, we have $(g_l(x))_{l\le d}=x$ for $x\in K$. So a randomized algorithm for $S$ on $F$ is
measurable in the sense of Section~\ref{sec:prelim} if and only if its points and its output are measurable
functions of $\omega$ and the observed values, and measurability is preserved in both directions, since
$a\mapsto x_a$ and $p^*$ are Borel measurable. Finally, every linear functional on the span of $F$ is of the form
$f_\xi\mapsto\ip a\xi$, so that $\Lambda^{\rm all}$ gives the same algorithms.

\smallskip
\emph{Bernstein widths.} Let $V\subset B(K)$ be a subspace of dimension $k+1$ and $\lambda>0$ with
$\lambda(B_{B(K)}\cap V)\subset F$. Then $V$ is contained in the span of $F$, that is, $V=\{f_\xi\colon\xi\in U\}$
for a subspace $U\subset\R^d$ of dimension $k+1$, and the inclusion says that $\xi\in B$ for all $\xi\in U$ with
$p(\xi)\le\lambda$. Conversely, every such $U$ gives a subspace $V$. If $B=B_p$, the inclusion holds for $\lambda=1$ and
fails for $\lambda>1$.
\end{proof}

So every example of a norm $p$ and a set $B$ for which randomization helps for the estimation of $p$ on $B$ from
linear measurements is an example for the problem of this paper. The classes of Section~\ref{sec:example} are of this
form with $B=B_1^d$ and the norms $\max\{\norm\cdot_\infty,\tau\norm\cdot_1\}$, $\norm\cdot_\infty$ and
$\norm\cdot_1$, and the problem of estimating a norm on its own unit ball, which is studied in
Section~\ref{sec:ran}, is the case $B=B_p$.

\subsection{Linear algorithms and measurability}
For linear problems with target space $B(D)$ on convex and symmetric classes, linear algorithms are optimal for
every fixed choice of the measurements. This goes back to Smolyak and Bakhvalov~\cite{Ba71} for real-valued
problems; for the general case see Creutzig and Wojtaszczyk~\cite{CW04}. We use it in the following form, which is stated
for function values in \cite[Section~1.3.1]{No88}, see also \cite[Proposition~1.2.6]{No88}.

\begin{lemma}\label{lem:linear}
Let $F\subset B(D)$ be nonempty, convex and symmetric, and let $x_1,\dots,x_n\in D$ be such that
$\rho:=\rho(x_1,\dots,x_n;F)$ is finite. Then there are functions $a_1,\dots,a_n\in B(D)$ with
\[
 \Bigl\|f-\sum_{i=1}^nf(x_i)\,a_i\Bigr\|_\infty\le\rho\qquad\text{for all }f\in F.
\]
\end{lemma}

A linear algorithm depends continuously on the function values. This gives algorithms that are measurable.

\begin{proposition}\label{prop:meas}
Let $F\subset B(D)$ be nonempty, convex and symmetric, $n\in\N$ and $T\in\{S,\Nf\}$. Then
$e_n^{\ran}(T,F)\le\rho_n(F)$. In particular,
\[
 e_n^{\ran}(S,F)\;\le\;2\,e_n^{\dt}(S,F).
\]
\end{proposition}

\begin{proof}
We may assume that $\rho_n(F)$ is finite. Let $x_1,\dots,x_n\in D$ be points for which the radius
$\rho:=\rho(x_1,\dots,x_n;F)$ is
finite, let $a_1,\dots,a_n$ be the functions from Lemma~\ref{lem:linear}, and put
\[
 g_y:=\sum_{i=1}^ny_i\,a_i\in B(D)\quad\text{for }y\in\R^n,\qquad L:=\max_{i\le n}\norm{a_i}_\infty .
\]
We consider the nonadaptive algorithm $A_n(f):=T(g_y)$ with $y=(f(x_1),\dots,f(x_n))$. Since
$\abs{T(f)-T(g)}\le\norm{f-g}_\infty$, its error is at most $\rho$ for every $f\in F$, and
\[
 \abs{T(g_y)-T(g_{y'})}\le\norm{g_y-g_{y'}}_\infty\le L\,\norm{y-y'}_1\qquad\text{for all }y,y'\in\R^n.
\]
So the output is a continuous function of the observed values. Since the points are fixed, $A_n$ is a measurable
randomized algorithm, which does not depend on $\omega$. Taking the infimum over the points gives
$e_n^{\ran}(T,F)\le\rho_n(F)$, and the last statement follows from Proposition~\ref{prop:det}.
\end{proof}

\section{Deterministic algorithms and Bernstein widths}\label{sec:det}

In this section we prove the deterministic part \eqref{eq:det-bern} of Theorem~\ref{thm:bernstein}. By
Proposition~\ref{prop:det} it is enough to bound the radius of information.

\begin{proposition}\label{prop:rho-bern}
Let $F\subset B(D)$ be nonempty, convex and symmetric. Then, for all $N,k\in\N$,
\[
 \rho_{N+k-1}(F)\;\le\;(N+k)\,\binom{N+k}{k}^{1/k}\;b_{N}(F).
\]
In particular, $\rho_{2N-1}(F)\le8N\,b_N(F)$ and $\rho_N(F)\le(N+1)^2\,b_N(F)$ for all $N\in\N$.
\end{proposition}

The two special cases are $k=N$, where we use $\binom{2N}N\le4^N$, and $k=1$. The bound \eqref{eq:det-bern}
follows immediately from the first one with $N=4n-1$ and Proposition~\ref{prop:det}. The second special case needs
no oversampling at all, at the price of a factor
$N^2$ instead of $N$.

Bounds of this kind have a long history for arbitrary linear measurements in place of function values, where the
radius of information is replaced by the Gelfand or Kolmogorov widths. Mityagin and Henkin~\cite{MH63} proved a
bound with the factor $(N+1)^2$ and conjectured that the factor $N+1$ is enough, and this was recently proved up to
the constant $e$ by Ullrich~\cite{Ul26}. The case $k=1$ of Proposition~\ref{prop:rho-bern} is the analogue of the
bound of Mityagin and Henkin for function values. Up to the factor $e$ it also follows from \cite{Ul26} together
with Novak's bound \cite[Proposition~1.2.5]{No88}, which says that in $B(D)$ the error of linear algorithms using
$n$ function values exceeds the Kolmogorov width by a factor of at most $n+1$; see also \cite{KPUU}. The bound with
oversampling seems to be new for a single index; with geometric means, a bound of this type is contained in
\cite[Theorem~1.2]{KNU}. Combining \cite{Ul26} with the sampling projections of \cite{KPUU}, which use $2N$
points and have norm of order $\sqrt N$, one only obtains a factor of order $N^{3/2}$ instead of $N$.

The proof of Proposition~\ref{prop:rho-bern} uses determinants of function values, namely, for bounded $F$ and
$j\in\N$, the numbers
\[
 \Gamma_j(F):=\sup\bigl\{\abs{\det[g_l(x_i)]_{i,l=1}^j}\colon g_1,\dots,g_j\in F,\ x_1,\dots,x_j\in D\bigr\},
\]
which are finite since $F$ is bounded. Determinants of this kind are a
classical tool in the theory of $s$-numbers, where they lead to bounds with geometric means. For function values,
the numbers $\Gamma_j(F)$ are, up to the normalization, the variant of the Grothendieck numbers that is
used in \cite[Section~5.3]{KU26}, see also \cite{KNU}. There, the minimal error of algorithms using function values
is bounded by these numbers with a triangular matrix of function values \cite[Lemma~5.14]{KU26}, and these numbers
are bounded by a geometric mean of the Hilbert numbers of $F$ \cite[Theorem~5.15]{KU26}, which gives
\cite[Theorem~1.2]{KNU}. Step~1 and part~(a) of Lemma~\ref{lem:det} below are variants of these two arguments.
A bound for a single index instead of a geometric mean was obtained in \cite{Ul26}, for arbitrary linear
measurements, by comparing two consecutive determinant quantities, and we follow this idea. The proof of
\cite{Ul26} uses linear combinations of the measurements and is not available for function values. What is new here is
that $\Gamma_{N+k}(F)$ and $\Gamma_N(F)$ are compared through the singular values of a matrix of function values with
almost maximal determinant, see Steps~3 and~4 in the proof of Lemma~\ref{lem:det}(b). This makes the comparison
possible for function values. For the
argument the class has to be bounded; the general case then follows by a truncation. For a square matrix $M$ we
write $s_1(M)\ge s_2(M)\ge\dots$ for its singular values.

\begin{lemma}\label{lem:det}
Let $F\subset B(D)$ be nonempty, bounded, convex and symmetric, let $N,k\in\N$, and let $g_1,\dots,g_{N+k}\in F$ and
$x_1,\dots,x_{N+k}\in D$ be such that the matrix $M:=[g_j(x_i)]_{i,j=1}^{N+k}$ is nonsingular.
\begin{enumerate}
\item[(a)] We have $b_{N}(F)\ge s_{N+1}(M)/(N+k)$.
\item[(b)] If $\abs{\det M}\ge(1-\eps)\,\Gamma_{N+k}(F)$ for some $0<\eps<1$, then
\[
 s_{N+1}(M)^k\;\ge\;(1-\eps)\binom{N+k}{k}^{-1}\rho_{N+k-1}(F)^k .
\]
\end{enumerate}
\end{lemma}

\begin{proof}
(a) Let $v_1,\dots,v_{N+k}$ be orthonormal right singular vectors of $M$, so that $\abs{Mv_i}=s_i(M)$ and
$Mv_i\perp Mv_j$ for $i\ne j$, and let $W:=\spn\{v_1,\dots,v_{N+1}\}$. For $w\in\R^{N+k}$ put
$f_w:=\sum_{j=1}^{N+k}w_jg_j$, so that $(f_w(x_i))_{i\le N+k}=Mw$. For $w\in W$ we have
$\abs{Mw}\ge s_{N+1}(M)\abs w$ and therefore
\[
 \norm{f_w}_\infty\ge\max_{i\le N+k}\abs{f_w(x_i)}=\norm{Mw}_\infty\ge\frac{\abs{Mw}}{\sqrt{N+k}}
 \ge\frac{s_{N+1}(M)\,\abs{w}}{\sqrt{N+k}}\ge\frac{s_{N+1}(M)}{N+k}\,\norm{w}_1.
\]
In particular $w\mapsto f_w$ is injective on $W$, and $V:=\{f_w\colon w\in W\}$ is an $(N+1)$-dimensional subspace
of $B(D)$. If $f=f_w\in V$ satisfies $\norm f_\infty\le s_{N+1}(M)/(N+k)$, then $\norm w_1\le1$, and hence
$f=\sum_jw_jg_j\in\conv\{\pm g_1,\dots,\pm g_{N+k}\}\subset F$, because $F$ is convex and symmetric. Thus
$\frac{s_{N+1}(M)}{N+k}(B_{B(D)}\cap V)\subset F$.

\smallskip
(b) \emph{Step 1.} We claim that
\begin{equation}\label{eq:onestep}
 \Gamma_j(F)\ge\rho_{j-1}(F)\,\Gamma_{j-1}(F)\qquad\text{for all }j\ge2.
\end{equation}
Let $0<\eta<1$ and choose $h_1,\dots,h_{j-1}\in F$ and $z_1,\dots,z_{j-1}\in D$ such that the matrix
$H:=[h_l(z_i)]_{i,l<j}$ satisfies $\abs{\det H}\ge(1-\eta)\Gamma_{j-1}(F)$. By the definition of $\rho_{j-1}(F)$,
there is $f\in F$ with $f(z_1)=\dots=f(z_{j-1})=0$ and
$\norm f_\infty\ge(1-\eta)\rho_{j-1}(F)$, and there is $z\in D$ with $\abs{f(z)}\ge(1-\eta)\norm f_\infty$. Let $H'$ be
the $j\times j$ matrix of the functions $h_1,\dots,h_{j-1},f$ evaluated at $z_1,\dots,z_{j-1},z$. Its last column
is $(0,\dots,0,f(z))^\top$, so $\det H'=f(z)\det H$, and therefore
$\Gamma_j(F)\ge\abs{\det H'}\ge(1-\eta)^3\rho_{j-1}(F)\Gamma_{j-1}(F)$. Letting $\eta\to0$ proves
\eqref{eq:onestep}. This is the argument of \cite[Lemma~5.14]{KU26}.

\smallskip
\emph{Step 2.} Applying \eqref{eq:onestep} for $j=N+1,\dots,N+k$ and using $\rho_{j-1}(F)\ge\rho_{N+k-1}(F)$ for
$j\le N+k$, we get
\begin{equation}\label{eq:telescope}
 \Gamma_{N+k}(F)\ge\Gamma_N(F)\prod_{j=N+1}^{N+k}\rho_{j-1}(F)\ge\Gamma_N(F)\,\rho_{N+k-1}(F)^k.
\end{equation}

\smallskip
\emph{Step 3.} We show that $\prod_{i\le N}s_i(M)\le\binom{N+k}{N}\,\Gamma_N(F)$. For $I,J\subset\{1,\dots,N+k\}$
with $\abs I=\abs J=N$ let $M_{I,J}$ be the submatrix of $M$ with rows $I$ and columns $J$. Let $Q$ be the
$(N+k)\times N$ matrix whose columns are the singular vectors $v_1,\dots,v_N$ from part (a), so that
$Q^\top Q=I_N$ and $Q^\top M^\top MQ$ is the diagonal matrix with entries $s_1(M)^2,\dots,s_N(M)^2$, and let $M_{I,:}$ and $Q_{J,:}$ be the submatrices of $M$ and $Q$
with rows $I$ and $J$, respectively. By the Cauchy--Binet formula, applied twice,
\[
 \det\bigl(Q^\top M^\top MQ\bigr)=\sum_{\abs I=N}\det\bigl(M_{I,:}Q\bigr)^2,
 \qquad
 \det\bigl(M_{I,:}Q\bigr)=\sum_{\abs J=N}\det M_{I,J}\,\det Q_{J,:}.
\]
Each $M_{I,J}$ is the matrix of the functions $g_j$, $j\in J$, at the points $x_i$, $i\in I$, hence
$\abs{\det M_{I,J}}\le\Gamma_N(F)$. By the Cauchy--Schwarz inequality and
$\sum_J\det(Q_{J,:})^2=\det(Q^\top Q)=1$ (Cauchy--Binet once more),
\[
 \abs{\det(M_{I,:}Q)}\le\Gamma_N(F)\sum_{\abs J=N}\abs{\det Q_{J,:}}\le\Gamma_N(F)\binom{N+k}{N}^{1/2},
\]
and hence $\prod_{i\le N}s_i(M)^2=\det(Q^\top M^\top MQ)\le\binom{N+k}N^2\Gamma_N(F)^2$. The inequality of this
step is classical: $\prod_{i\le N}s_i(M)$ is the largest singular value of the $N$-th compound matrix
$[\det M_{I,J}]_{I,J}$ of $M$, which has $\binom{N+k}N$ rows and columns and entries of absolute value at most
$\Gamma_N(F)$.

\smallskip
\emph{Step 4.} Note that $\Gamma_N(F)>0$, since the nonsingular matrix $M$ has a nonsingular $N\times N$
submatrix. By Step 3, the assumption on $M$, and \eqref{eq:telescope},
\begin{align*}
 s_{N+1}(M)^k&\ge\prod_{i=N+1}^{N+k}s_i(M)=\frac{\abs{\det M}}{\prod_{i\le N}s_i(M)}\\
 &\ge\frac{(1-\eps)\,\Gamma_{N+k}(F)}{\binom{N+k}N\,\Gamma_N(F)}\ge(1-\eps)\binom{N+k}{k}^{-1}\rho_{N+k-1}(F)^k. \qedhere
\end{align*}
\end{proof}

\begin{proof}[Proof of Proposition~\ref{prop:rho-bern}]
We first assume that $F$ is bounded. If $F$ does not contain $N+k$ linearly independent functions, then $\spn F$ has dimension at most $N+k-1$, and
there are $N+k-1$ points in $D$ such that every function in $\spn F$ that vanishes at these points is zero;
hence $\rho_{N+k-1}(F)=0$. Otherwise there are
$g_1,\dots,g_{N+k}\in F$ and points with a nonsingular matrix of function values (for linearly independent
functions $f_1,\dots,f_r$ the vectors $(f_1(x),\dots,f_r(x))$, $x\in D$, span $\R^r$), so $\Gamma_{N+k}(F)>0$. For
$0<\eps<1$ choose $g_j$ and $x_i$ with $\abs{\det M}\ge(1-\eps)\Gamma_{N+k}(F)$. By Lemma~\ref{lem:det},
\[
 \rho_{N+k-1}(F)\le(1-\eps)^{-1/k}\binom{N+k}{k}^{1/k}s_{N+1}(M)\le(1-\eps)^{-1/k}\binom{N+k}{k}^{1/k}(N+k)\,b_{N}(F),
\]
and $\eps\to0$ proves the claim for bounded classes.

If $F$ is unbounded, we consider the classes $F_R:=\{f\in F\colon\norm f_\infty\le R\}$ for $R>0$, which are
nonempty, bounded, convex and symmetric. If $h\in F$ vanishes at given points $x_1,\dots,x_n$, then the functions
$sh$ with $0\le s\le1$ belong to $F$, since $0\in F$ and $F$ is convex, and they vanish at these points. Hence
$\rho(x_1,\dots,x_n;F_R)=\min\{\rho(x_1,\dots,x_n;F),R\}$ and therefore $\rho_n(F_R)=\min\{\rho_n(F),R\}$ for all
$n$. Moreover, $b_N(F_R)\le b_N(F)$ since $F_R\subset F$. The bound for the bounded class $F_R$ gives
\[
 \min\bigl\{\rho_{N+k-1}(F),R\bigr\}\le(N+k)\,\binom{N+k}{k}^{1/k}\,b_N(F)\qquad\text{for all }R>0,
\]
and $R\to\infty$ proves the claim.
\end{proof}

\begin{remark}\label{rem:sharp}
The factor $N$ in the bound $\rho_{2N-1}(F)\le8N\,b_N(F)$ cannot be improved in order. Let $D=\{1,\dots,d\}$
and let $F$ be the unit ball of $\ell_1^d$, considered as a class of functions on $D$. Then $\rho_{m}(F)=1$ for
$m<d$, since a unit vector at a coordinate that is not sampled belongs to $F$. On the other hand,
\[
 b_{N}(F)=\frac1{N+1}\qquad\text{for }N<d .
\]
The lower bound is obtained for the space $V$ of all vectors supported in the first $N+1$ coordinates. For the upper
bound let $V\subset\R^d$ be a subspace of dimension $N+1$ and $\lambda>0$ with
$\lambda(B_\infty^d\cap V)\subset B_1^d$. The polytope $B_\infty^d\cap V$ is bounded and hence has a vertex $w$. At
least $N+1$ of the constraints $\abs{w_i}\le1$ are active at $w$, since otherwise there is a nonzero $u\in V$ with
$u_i=0$ for all active $i$, and $w\pm tu$ belongs to the polytope for small $t>0$. Hence $\norm w_1\ge N+1$, and
$\lambda w\in B_1^d$ implies $\lambda\le1/(N+1)$.
So $\rho_{2N-1}(F)=(N+1)\,b_N(F)$ for $d\ge2N$.

Likewise, the factor $n$ in \eqref{eq:det-bern} cannot be
improved: for $d\ge8n$ we have $e^{\dt}_{8n-1}(S,F)\ge\frac12\rho_{8n-1}(F)=\frac12=2n\,b_{4n-1}(F)$. We do not
know whether the factor $(N+1)^2$ in the bound $\rho_N(F)\le(N+1)^2b_N(F)$ can be replaced by a constant multiple
of $N$; in the present example $\rho_N(F)=(N+1)\,b_N(F)$. This is the analogue for function values of the
conjecture of Mityagin and Henkin~\cite{MH63}, and it does not follow by combining \cite{Ul26} with Novak's bound, since the factor $n+1$
in the latter cannot be improved, see \cite[Proposition~1.2.5]{No88} and \cite{KPUU}. The unit ball of $\ell_1^d$ is also the extremal example
in \cite{Ul26}.
\end{remark}

\begin{remark}
The proof shows more than stated. In part (a) one can measure the ball in $L_2$ of the uniform
distribution on $x_1,\dots,x_{N+k}$ instead of the uniform norm, which gives the bound with the Hilbert numbers of
$F$ in place of the Bernstein widths; this is the construction in the proof of \cite[Theorem~5.15]{KU26}, see
\cite{KNU,Ul26} for these notions. We do not need this here.
\end{remark}

\section{Randomized algorithms and Bernstein widths}\label{sec:ran}

In this section we prove the lower bound \eqref{eq:ran-bern} for randomized algorithms. If $F$ contains a ball of
dimension $2n$, then the problem of computing $\norm f_\infty$ on this ball is the problem of estimating a norm on
its own unit ball in $\R^{2n}$, where function values are linear measurements. We put a prior on this ball and
show that $n$ adaptive linear measurements leave an uncertainty at least of order $1/n$ about the norm.

\subsection{Linear measurements and fibres}
Let $d\in\N$ and $k<d$. A (deterministic, adaptive) algorithm using $k$ linear measurements of an input $x\in\R^d$
is given by measurable mappings $q_j\colon\R^{j-1}\to\R^d$, $j\le k$, and $\varphi\colon\R^k\to\R$. For the input $x$
it computes
\[
 y_j=\ip{q_j(y_1,\dots,y_{j-1})}{x},\quad j=1,\dots,k,
\]
and returns $A(x)=\varphi(y_1,\dots,y_k)$. We may
and do assume that the measurement vectors $q_1$, $q_2(y_1)$, \dots, $q_k(y_1,\dots,y_{k-1})$ are orthonormal for
all $y$. Indeed, we
can apply the Gram--Schmidt procedure along each path of observed values, where a vector that is linearly
dependent on the previous ones is replaced by a unit vector that is orthogonal to them and chosen by a fixed
rule (this is possible since $k<d$). The values observed by the original algorithm are measurable functions of the values observed by the new
one, so that the new algorithm can simulate the original one, see also \cite{KNU2,TWW}. For $y\in\R^k$ we put
\[
 a_y:=\sum_{j\le k}y_j\,q_j(y_1,\dots,y_{j-1}),\qquad W_y:=\{q_1,q_2(y_1),\dots,q_k(y_1,\dots,y_{k-1})\}^\perp ,
\]
so that $W_y$ is a subspace of dimension $d-k$.

\begin{lemma}\label{lem:fibres}
The inputs that produce the values $y\in\R^k$ are exactly the points of the affine subspace $a_y+W_y$, and for every
nonnegative measurable function $h$ on $\R^d$,
\begin{equation}\label{eq:fubini}
 \int_{\R^d}h(x)\,dx=\int_{\R^k}\int_{W_y}h(a_y+w)\,dw\,dy,
\end{equation}
where $dw$ is the $(d-k)$-dimensional Lebesgue measure on $W_y$. Consequently, if the input $x$ is random with
density $\gamma$, then the vector $y$ of observed values has the density $y\mapsto\int_{W_y}\gamma(a_y+w)\,dw$, and conditionally
on $y$ the input is distributed on $a_y+W_y$ with density proportional to $w\mapsto\gamma(a_y+w)$.
(All functions of $y$ that occur here are measurable, since the inner integral can be written as an integral over
$\R^d$ of a jointly measurable function with the help of the orthogonal projection onto $W_y$.)
\end{lemma}

\begin{proof}
If $x$ produces $y$, then $\ip{q_j(y_1,\dots,y_{j-1})}{x}=y_j$ for all $j$, that is, $x-a_y\in W_y$, and vice versa.
The formula \eqref{eq:fubini} is proved by induction on $k$. For $k=1$ it is Fubini's theorem in an orthonormal basis
containing $q_1$. For $k>1$ we write $x=y_1q_1+x'$ with $x'\in q_1^\perp$. For fixed $y_1$, the vectors
$q_2(y_1),\dots,q_k(y_1,\cdot)$ define an algorithm with $k-1$ orthonormal measurements on the $(d-1)$-dimensional
space $q_1^\perp$, and we apply the induction hypothesis to $x'\mapsto h(y_1q_1+x')$ and integrate over $y_1$.
\end{proof}

\subsection{Estimating the gauge of a convex body}
For a symmetric convex body $C\subset\R^{2n}$ (compact, convex, $C=-C$, with nonempty interior) let
$\norm{\xi}_C:=\inf\{s>0\colon\xi\in sC\}$ be its gauge, a norm on $\R^{2n}$ with unit ball $C$.

\begin{lemma}[Gauge estimation]\label{lem:gauge}
Let $n\ge100$ and let $C\subset\R^{2n}$ be a symmetric convex body. Let $u$ be uniformly distributed on $C$, let $t$
be uniformly distributed on $[\frac12,1]$ and independent of $u$, and put $\xi:=tu$ and $\varrho:=\norm{\xi}_C$.
Then every algorithm $A$ using $n$ linear measurements satisfies
\[
 \E\abs{\varrho-A(\xi)}\ge\frac{c_0}{n},\qquad c_0=\frac1{175}.
\]
\end{lemma}

For linear problems, lower bounds for randomized algorithms in terms of Bernstein widths were obtained with the
uniform distribution on a ball and with Gaussian measures as priors, see Kunsch~\cite{Ku16,Ku17}. We do
not know whether a lower bound of order $1/n$ also holds for the uniform distribution on $C$, without the random
scale $t$, but this is irrelevant for our purpose, since we are free to choose the probability measure.

The idea of the proof is as follows. Given the observed values $y$, the input lies in the $n$-dimensional affine
subspace $a_y+W_y$, and we compare the conditional probability that $\varrho$ lies in an interval $J=[m-\eps,m+\eps]$ around a possible
output $m=\varphi(y)$ of the algorithm, where $\eps$ is of order $1/n$, with that of an interval $J'$ of length of order $1/n$ directly above it. The volume of the inputs
with $\varrho\le s$ in this subspace grows at least like $s^n$, by the Brunn--Minkowski inequality, so $J'$ contains
at least as much volume as everything below it. The density of $\xi$ decays like $\varrho^{1-2n}$, so on $J'$ it is
smaller than on $J$ by at most a constant factor. Hence $J'$ carries a fixed fraction of the mass of $J$, and
$\varrho$ lies outside of $J$ with a probability that is bounded from below. So no matter how we choose
$\varphi(y)$, we miss $\varrho$ by more than $\eps$ with constant probability. This argument needs room above $J$ and
therefore fails if $m$ is within order $1/n$ of the maximal value one; but because of the random scale, $\varrho$ is
rarely that large, so that such an output also misses $\varrho$ for most $y$. We first collect the facts about the prior and the
posterior.

\begin{lemma}\label{lem:prior}
In the setting of Lemma~\ref{lem:gauge}, define
\[
 \Pi(r):=\int_{\max\{r,1/2\}}^1t^{-2n}\,dt\quad(0\le r\le1),\qquad\Pi(r):=0\quad(r>1),
\]
and, for fixed $y\in\R^n$, let $L_y:=a_y+W_y$ be the $n$-dimensional affine subspace of all inputs that produce the
values $y$, see Lemma~\ref{lem:fibres}, and let $\Theta(s):=\vol_n(sC\cap L_y)$ for $s\ge0$. Then:
\begin{enumerate}
\item[(i)] The random vector $\xi$ has the density $z\mapsto\frac{2}{\vol(C)}\,\Pi(\norm z_C)$ on $\R^{2n}$.
\item[(ii)] For almost every $y$ with respect to the distribution $P$ of the observed values, the conditional
distribution $Q_y$ of $\varrho$ given $y$ satisfies
\[
 Q_y([0,s])=\frac{\int_{[0,s]}\Pi\,d\Theta}{\int_{[0,1]}\Pi\,d\Theta}\qquad(0\le s\le1),
\]
where $d\Theta$ denotes the measure on $[0,\infty)$ with distribution function $\Theta$, that is, the measure
$\nu$ with $\nu([0,s])=\Theta(s)$ for all $s\ge0$.
\item[(iii)] The function $s\mapsto\Theta(s)/s^n$ is nondecreasing on $(0,\infty)$. In particular,
$\Theta(s')\le(s'/s)^n\,\Theta(s)$ for $0<s'\le s$.
\item[(iv)] $\Pi$ is continuous and nonincreasing, and for $n\ge100$, $\beta>0$ and $0\le x\le x'$ with
$x'\le e^{\beta/n}x$ and $x'(1+\frac2n)\le1$,
\[
 \frac{\Pi(x)}{\Pi(x')}\le\frac{50\,e^{2\beta}-1}{49}=:R_\beta .
\]
\end{enumerate}
\end{lemma}

\begin{proof}
(i) Given $t$, the vector $\xi=tu$ is uniformly distributed on $tC$, with density
$z\mapsto t^{-2n}\mathbf 1[\norm z_C\le t]/\vol(C)$. Integrating against the density $2\cdot\mathbf 1_{[1/2,1]}(t)$ of $t$
gives the claim.

\smallskip
(ii) By (i) and Lemma~\ref{lem:fibres}, the conditional distribution of $\xi$ given $y$ has a density on $L_y$ that
is proportional to $z\mapsto\Pi(\norm z_C)$. The normalizing constant $\int_{L_y}\Pi(\norm z_C)\,dz$ is finite, since
$\Pi$ is bounded and vanishes outside of the bounded set $C$, and it is positive for $P$-almost every $y$, since
it is a constant multiple of the density of $P$ at $y$. Let $\nu$ be the image of the $n$-dimensional Lebesgue measure on $L_y$ under the
mapping $z\mapsto\norm z_C$, that is, $\nu(B):=\vol_n\{z\in L_y\colon\norm z_C\in B\}$ for Borel sets
$B\subset[0,\infty)$. Since $\{z\colon\norm z_C\le s\}=sC$, we have $\nu([0,s])=\vol_n(sC\cap L_y)=\Theta(s)$. So
$\nu$ is the measure $d\Theta$ of the statement, and in particular $\nu((s,s'])=\Theta(s')-\Theta(s)$ for $s<s'$.
By the change of variables formula for image measures,
\[
 \int_{L_y}h(\norm z_C)\,dz=\int_{[0,\infty)}h\,d\nu
\]
for every nonnegative measurable function $h$ on $[0,\infty)$. For $h=\Pi\cdot\mathbf 1_{[0,s]}$ this gives
\[
 \int_{\{z\in L_y\colon\norm z_C\le s\}}\Pi(\norm z_C)\,dz=\int_{[0,s]}\Pi\,d\Theta .
\]
The left-hand side, divided by the normalizing constant, is $Q_y([0,s])$, and $s=1$ gives the normalizing
constant itself, since $\Pi$ vanishes on $(1,\infty)$.

\smallskip
(iii) We have $sC\cap(a_y+W_y)=s\bigl(C\cap(s^{-1}a_y+W_y)\bigr)$, hence $\Theta(s)/s^n=\psi(s^{-1})$ with
$\psi(\lambda):=\vol_n\bigl(C\cap(\lambda a_y+W_y)\bigr)$. The function $\psi$ is even, since $C=-C$, and
$\psi^{1/n}$ is concave on the closed interval of all $\lambda$ for which the section
$C\cap(\lambda a_y+W_y)$ is nonempty, by the Brunn--Minkowski inequality in the
$n$-dimensional space $W_y$ and the convexity of $C$, see \cite{Ga02}. A nonnegative even function that is concave
on a symmetric interval and zero outside is nonincreasing in $\abs\lambda$. Hence $\psi(s^{-1})$ is nondecreasing
in $s$.

\smallskip
(iv) Continuity and monotonicity are clear. If $x'\le\frac12$, the ratio is $1\le R_\beta$. If $x<\frac12<x'$, then
$\Pi(x)=\Pi(\frac12)$ and we may replace $x$ by $\frac12$. So let $\frac12\le x\le x'$. For $s\in[\frac12,1]$ we have
$\Pi(s)=(s^{1-2n}-1)/(2n-1)$, and hence
\[
 \frac{\Pi(x)}{\Pi(x')}=\frac{x^{1-2n}-1}{(x')^{1-2n}-1}=\frac{A-u}{1-u}
 \qquad\text{with}\quad A:=\Bigl(\frac{x'}x\Bigr)^{2n-1},\quad u:=(x')^{2n-1}.
\]
Here $1\le A\le e^{2\beta}$. Using $\log(1+s)\ge s-\frac{s^2}2$ and $n\ge100$, we get
$(2n-1)\log(1+\frac2n)\ge(2n-1)(\frac2n-\frac2{n^2})\ge4-\frac6n\ge3.94$, so that
$u\le(1+\frac2n)^{-(2n-1)}\le e^{-3.94}\le\frac1{50}$. Since $u\mapsto(A-u)/(1-u)$ is nondecreasing on $[0,1)$
for $A\ge1$, the ratio is at most $(e^{2\beta}-\frac1{50})/(1-\frac1{50})=R_\beta$.
\end{proof}

\begin{proof}[Proof of Lemma~\ref{lem:gauge}]
Let $\varphi(y)$ be the output of the algorithm and let $P$ be the distribution of $y$. Then
\[
 \E\abs{\varrho-A(\xi)}=\int\E\bigl[\abs{\varrho-\varphi(y)}\,\big|\,y\bigr]\,dP(y)
 \ge\int_Y\inf_{m\in\R}\int\abs{s-m}\,dQ_y(s)\,dP(y)
\]
for every measurable set $Y\subset\R^n$. Here $(Q_y)$ is a regular conditional distribution of $\varrho$ given $y$,
so that $y\mapsto Q_y(B)$ is measurable for Borel sets $B$, and the infimum may be restricted to rational $m$. In
the following we only consider those $y$ for which the formula of Lemma~\ref{lem:prior}(ii) holds; this is a set
of full measure with respect to $P$. We choose $Y:=Y_1\cap Y_2$ with
\[
 Y_1:=\bigl\{y\colon Q_y\bigl([0,1-\tfrac5n]\bigr)\ge\tfrac12\bigr\},\qquad
 Y_2:=\bigl\{y\colon Q_y\bigl([\tfrac9{20},1]\bigr)\ge\tfrac12\bigr\}.
\]

\smallskip
\emph{Step 1: $P(Y)\ge0.79$.} Since $\varrho\le t$, we have $\P(\varrho>1-\frac5n)\le\P(t>1-\frac5n)=\frac{10}n$. By
Markov's inequality, applied to the function $y\mapsto Q_y((1-\frac5n,1])$ whose mean is $\P(\varrho>1-\frac5n)$, we get
$P(Y_1^c)\le\frac{20}n\le0.2$. Since $\varrho\ge\frac12\norm u_C$ and
$\P(\norm u_C<\frac9{10})=\vol(\frac9{10}C)/\vol(C)=(\frac9{10})^{2n}\le10^{-9}$, the same argument gives
$P(Y_2^c)\le2\cdot10^{-9}$. Hence $P(Y)\ge1-0.2-2\cdot10^{-9}\ge0.79$.

\smallskip
\emph{Step 2: anti-concentration for $\frac25\le m\le1-\frac4n$.} Fix $y$ and such an $m$, and put
$\eps:=\frac m{4n}$, $l:=m-\eps$, $r:=m+\eps$ and $J:=[l,r]$. We show that
\begin{equation}\label{eq:anti}
 Q_y(J)\le1-0.074.
\end{equation}
Let $\mu$ be the measure $\Pi\,d\Theta$ on $[0,1]$, so that $Q_y=\mu/\mu([0,1])$. Since $\Pi$ is nonincreasing, we
have $\mu(J)\le\Pi(l)\,\Theta(r)$. If $\Theta(r)=0$, then $Q_y(J)=0$ and \eqref{eq:anti} holds. Otherwise let
$r':=r(1+\frac1n)$ and $J':=(r,r']$.
By Lemma~\ref{lem:prior}(iii) and $(1+\frac1n)^n\ge e^{1-1/(2n)}\ge2.7$, we have $\Theta(r')\ge2.7\,\Theta(r)$, hence
\[
 \mu(J')\ge\Pi(r')\bigl(\Theta(r')-\Theta(r)\bigr)\ge1.7\;\Pi(r')\,\Theta(r).
\]
We have $r'/l=\frac{1+1/(4n)}{1-1/(4n)}(1+\frac1n)\le e^{1.51/n}$, because
$\log\frac{1+x}{1-x}\le2.04\,x$ for $0\le x\le\frac1{400}$ and $\log(1+\frac1n)\le\frac1n$, and
$r'(1+\frac2n)\le(1-\frac4n)(1+\frac1{4n})(1+\frac1n)(1+\frac2n)\le1$, since $1+x\le e^x$ for all $x\in\R$. Lemma~\ref{lem:prior}(iv) with
$\beta=1.51$ therefore gives $\Pi(l)\le R\,\Pi(r')$ with $R:=R_{1.51}=(50\,e^{3.02}-1)/49\le21$; in particular
$\Pi(r')>0$. With $w:=\Pi(r')\,\Theta(r)>0$ we thus have $\mu(J)\le R\,w$ and $\mu(J')\ge1.7\,w$. Since $J$ and $J'$
are disjoint subsets of $[0,1]$ and $a/(a+b)$ is increasing in $a$ and decreasing in $b$,
\[
 Q_y(J)\le\frac{\mu(J)}{\mu(J)+\mu(J')}\le\frac R{R+1.7}\le\frac{21}{22.7}\le1-0.074.
\]
This proves \eqref{eq:anti}, and consequently
$\int\abs{s-m}\,dQ_y(s)\ge0.074\,\eps=0.074\cdot\frac m{4n}\ge\frac{0.0074}n$.

\smallskip
\emph{Step 3: the remaining values of $m$.} For $y\in Y$ and $m<\frac25$ we have
$\int\abs{s-m}\,dQ_y(s)\ge Q_y([\frac9{20},1])\cdot\frac1{20}\ge\frac1{40}$. For $m>1-\frac4n$ we have
$\int\abs{s-m}\,dQ_y(s)\ge Q_y([0,1-\frac5n])\cdot\frac1n\ge\frac1{2n}$.

\smallskip
\emph{Conclusion.} For $y\in Y$ we have shown $\inf_m\int\abs{s-m}\,dQ_y(s)\ge\frac{0.0074}n$, and with Step~1 we
obtain $\E\abs{\varrho-A(\xi)}\ge0.79\cdot\frac{0.0074}n\ge\frac1{175\,n}$.
\end{proof}

\begin{remark}
The order $1/n$ in Lemma~\ref{lem:gauge} cannot be improved for this prior, since there are bodies and algorithms
with an error of this order. Let $C=B_2^n\times B_\infty^n$ and write $\xi=(a,b)$ with $a,b\in\R^n$. The
algorithm that measures the $n$ coordinates of $b$ and returns $\norm b_\infty$ has the error
$\varrho-\norm b_\infty=(\abs a-\norm b_\infty)_+$. For $u$ uniformly distributed on $C$, the Euclidean norm of
its first component and the maximum norm of its second component are independent, each with the distribution
function $x\mapsto x^n$ on $[0,1]$. The maximum of two such variables has the mean $\frac{2n}{2n+1}$, and each of
them has the mean $\frac n{n+1}$. Since $\E[t]=\frac34$, this gives
\[
 \E\,\abs{\varrho-\norm b_\infty}=\frac34\Bigl(\frac{2n}{2n+1}-\frac n{n+1}\Bigr)=\frac34\cdot\frac{n}{(2n+1)(n+1)}\le\frac3{8n}.
\]
\end{remark}

\subsection{Proof of Theorem~\ref{thm:bernstein} and Theorem~\ref{thm:upper}}

\begin{lemma}\label{lem:ran-norm}
Let $F\subset B(D)$ be nonempty, convex and symmetric, and let $n\ge100$. Then
\[
 e_n^{\ran}(\Nf,F)\ge e_n^{\ran}(\Nf,F,\Lambda^{\rm all})\ge\frac{c_0}n\,b_{2n-1}(F).
\]
\end{lemma}

\begin{proof}
We may assume that $b_{2n-1}(F)>0$. Let $0<\lambda<b_{2n-1}(F)$ and choose a subspace $V\subset B(D)$ with $\dim V=2n$ and
$\lambda(B_{B(D)}\cap V)\subset F$. Fix a basis $g_1,\dots,g_{2n}$ of $V$ and identify $\xi\in\R^{2n}$ with
$f_\xi:=\sum_j\xi_jg_j$. Then $C:=\{\xi\in\R^{2n}\colon\norm{f_\xi}_\infty\le\lambda\}$ is a symmetric convex body with
$\norm\xi_C=\norm{f_\xi}_\infty/\lambda$, and $f_\xi\in F$ for all $\xi\in C$. Let $(A^\omega)$ be a measurable
randomized algorithm for $\Nf$ on $F$ using $n$ linear measurements. Note that $V$ is spanned by functions
from $F$, and that $L(f_\xi)=\ip{(L(g_j))_j}{\xi}$ for every $L\in\Lambda^{\rm all}$. Hence, on the inputs $f_\xi$,
each $A^\omega$ is an algorithm using $n$ linear measurements of $\xi\in\R^{2n}$ in the sense of the beginning of this
section, with the measurement vectors $q_j^\omega(y):=\bigl(\psi_j^\omega(y)(g_l)\bigr)_{l\le2n}$ and the output
mapping $\varphi^\omega$; these are measurable for every fixed $\omega$, as sections of jointly measurable mappings. Let $\xi=tu$ be the random input of Lemma~\ref{lem:gauge}. By
\eqref{eq:bakhvalov} and Lemma~\ref{lem:gauge},
\[
 \sup_{f\in F}\E_\omega\abs{\Nf(f)-A^\omega(f)}\ge\E_\omega\E_\xi\bigl|\lambda\norm\xi_C-A^\omega(f_\xi)\bigr|
 =\lambda\,\E_\omega\E_\xi\Bigl|\norm\xi_C-\frac{A^\omega(f_\xi)}\lambda\Bigr|\ge\frac{\lambda\,c_0}n.
\]
Taking the infimum over all algorithms and letting $\lambda\to b_{2n-1}(F)$ proves the claim.
\end{proof}

\begin{proof}[Proof of Theorem~\ref{thm:bernstein}]
The bound \eqref{eq:det-bern} was proved after Proposition~\ref{prop:rho-bern}. For \eqref{eq:ran-bern} let
$n\ge50$. By Lemma~\ref{lem:StoN} and Lemma~\ref{lem:ran-norm}, applied with $2n\ge100$ measurements,
\[
 e_n^{\ran}(S,F)\ge e_n^{\ran}(S,F,\Lambda^{\rm all})\ge\tfrac12\,e_{2n}^{\ran}(\Nf,F,\Lambda^{\rm all})
 \ge\frac12\cdot\frac{c_0}{2n}\,b_{4n-1}(F). \qedhere
\]
\end{proof}

\begin{proof}[Proof of Theorem~\ref{thm:upper}]
By \eqref{eq:det-bern} and \eqref{eq:ran-bern},
\[
 e_{8n-1}^{\dt}(S,F)\le32n\,b_{4n-1}(F)\le32n\cdot\frac{4n}{c_0}\,e_n^{\ran}(S,F,\Lambda^{\rm all})
 \le\frac{128}{c_0}\,n^2\,e_n^{\ran}(S,F),
\]
and $128/c_0=22\,400$.
\end{proof}

\begin{remark}[Uniform approximation]\label{rem:app-upper}
The deterministic algorithm in Theorem~\ref{thm:upper} may even be asked to solve the harder problem of uniform
approximation. By Remark~\ref{rem:app}, Proposition~\ref{prop:rho-bern} and the bound \eqref{eq:ran-bern}, for all $n\ge50$,
\begin{equation}\label{eq:app-upper}
 e_{8n-1}^{\dt}(\APP,F,\Lambda^{\rm std})\;\le\;32\,n\,b_{4n-1}(F)\;\le\;22\,400\,n^2\,e_n^{\ran}(S,F,\Lambda^{\rm all}).
\end{equation}
Note that \eqref{eq:app-upper} compares two different problems, the hardest one for the deterministic
and the easiest one for the randomized algorithm. Moreover, the deterministic algorithm uses function values, while
the randomized algorithm may use arbitrary linear measurements. If uniform approximation is considered on both sides, then a
factor of order $n$ is enough. It is a linear problem, and for linear problems on convex classes the error of
randomized algorithms that use $n$ function values, or other continuous linear measurements, is at least
$\frac1{30}\,b_{2n-1}(F)$, see \cite[Theorem~9.23]{KU26}, which is based on \cite{Ku16,Ku17}. Together with
Remark~\ref{rem:app} and Proposition~\ref{prop:rho-bern} for $N=2n-1$ this shows that
\[
 e_{4n-1}^{\dt}(\APP,F,\Lambda^{\rm std})\;\le\;16\,n\,b_{2n-1}(F)\;\le\;480\,n\;e_n^{\ran}(\APP,F,\Lambda^{\rm lin}),
\]
where $e_n^{\ran}(\APP,F,\Lambda^{\rm lin})$ denotes the minimal error of randomized algorithms for uniform approximation
on $F$ that use $n$ continuous linear measurements, for example function values, as defined in \cite{KU26}; for
unbounded $F$ we apply the second inequality to the bounded
classes $F_R$ from the proof of Proposition~\ref{prop:rho-bern} and use that $F_R\subset F$ and
$b_{2n-1}(F_R)=\min\{b_{2n-1}(F),R\}$. A bound of this type, with a geometric mean of the
errors of randomized algorithms instead of a single-index comparison, is contained in \cite[Theorem~1.2]{KNU}.
\end{remark}

\section{A class with a gain of order \texorpdfstring{$n^{3/2}$}{n\textasciicircum(3/2)}}\label{sec:example}

In this section we prove Theorem~\ref{thm:lower}. We do this in the setting of Lemma~\ref{lem:dict}, that is, for
the estimation of a norm from linear measurements. For $d\in\N$ and $\tau\in(0,1]$ we consider the norm
\begin{equation}\label{eq:S-example}
 p_\tau(\xi):=\max\bigl\{\norm\xi_\infty,\ \tau\norm\xi_1\bigr\},\qquad\xi\in\R^d,
\end{equation}
on the set $B_1^d$, where it is bounded by one. We write $F_{d,\tau}$ for the corresponding class of
Lemma~\ref{lem:dict}, with $p=p_\tau$ and $B=B_1^d$. This class is convex, symmetric and bounded, and for all $n\in\N$,
\begin{equation}\label{eq:translate}
 e_n^{\dt}(S,F_{d,\tau})=e_n^{\dt}(p_\tau,B_1^d,\Lambda^{\rm all})\qquad\text{and}\qquad
 e_n^{\ran}(S,F_{d,\tau})=e_n^{\ran}(p_\tau,B_1^d,\Lambda^{\rm all}).
\end{equation}

For small $\tau$ the problem behaves like the estimation of $\norm\xi_\infty$ on $B_1^d$.
There, a randomized algorithm has to find the largest coordinate of $\xi$. Randomized algorithms that find the
large coordinates of a vector from few linear measurements are well known, in the literature on data streams, see
for example \cite{CM05}, and in information-based complexity, see \cite{KNW,KW}. They lead to an error of order
$\log(n)/n$, see Remark~\ref{rem:norms}(1). For $\tau=1$ we have
$p_1(\xi)=\norm\xi_1$, which can be estimated with an error of order $n^{-1/2}$ by a Monte Carlo method, see
Remark~\ref{rem:norms}(2). For deterministic algorithms the problem is hard for every $\tau$ if the dimension is large,
see Proposition~\ref{prop:ex-det}. The idea is to choose $\tau$ so small that the Monte Carlo error for the second
term in \eqref{eq:S-example}, which is of order $\tau n^{-1/2}$, is small, but so large that the first term only
matters if $\xi$ has a coordinate that a randomized algorithm can find and then measure \emph{exactly}.

\subsection{Deterministic algorithms}

\begin{proposition}\label{prop:ex-det}
Let $d>6^m$. Then for all $x_1,\dots,x_m\in\R^d$ there is $\xi\in B_1^d$ with $\ip{x_i}\xi=0$ for all $i\le m$ and
$\norm\xi_\infty\ge\frac14$. Consequently, $e_m^{\dt}(\norm\cdot_\infty,B_1^d,\Lambda^{\rm all})\ge\frac18$.
\end{proposition}

The second statement follows from the first: a deterministic algorithm that uses the measurements $\ip{x_i}\cdot$
for the input zero cannot distinguish such a vector $\xi$ from zero. It returns the same number for both inputs and
therefore has an error of at least $\frac12\norm\xi_\infty\ge\frac18$. The same argument applies to every norm $p$ on
$\R^d$ with $p\ge\norm\cdot_\infty$. In particular,
\begin{equation}\label{eq:ex-det-tau}
 e_m^{\dt}(p_\tau,B_1^d,\Lambda^{\rm all})\ge\frac18\qquad\text{for all }\tau\in(0,1]\text{ and }d>6^m.
\end{equation}
The proof of the first statement is a packing argument. We need the following duality.

\begin{lemma}\label{lem:duality}
Let $V\subset\R^d$ be a linear subspace and $e\in\R^d$. Then there is $\xi\in V^\perp$ with $\norm\xi_1\le1$ and
$\ip e\xi=\dist_\infty(e,V):=\inf_{v\in V}\norm{e-v}_\infty$.
\end{lemma}

\begin{proof}
Let $\delta:=\dist_\infty(e,V)$. If $\delta=0$, take $\xi=0$. Otherwise the convex sets $V$ and
$\{e+\delta b\colon\norm b_\infty<1\}$ are disjoint, and by the separation theorem there is $\xi\ne0$ with
$\ip\xi v\le\ip\xi{e+\delta b}$ for all $v\in V$ and $\norm b_\infty<1$. Since $V$ is a subspace, $\xi\in V^\perp$,
and then $0\le\ip\xi e+\delta\ip\xi b$ for all $\norm b_\infty<1$ gives $\ip\xi e\ge\delta\norm\xi_1$. After
normalizing $\norm\xi_1=1$ we have $\ip e\xi\ge\delta$, and a multiple of $\xi$ gives equality.
\end{proof}

\begin{proof}[Proof of Proposition~\ref{prop:ex-det}]
Let $V:=\spn\{x_1,\dots,x_m\}$ and $k:=\dim V\le m$. Suppose that
$\dist_\infty(e_j,V)<\frac14$ for every unit vector $e_j$, $j\le d$, and choose $v_j\in V$ with
$\norm{e_j-v_j}_\infty<\frac14$. Then $\norm{v_j}_\infty<\frac54$ and
$\norm{v_j-v_{j'}}_\infty>\norm{e_j-e_{j'}}_\infty-\frac12=\frac12$ for $j\ne j'$. Thus the sets
$v_j+\frac14(B_\infty^d\cap V)$, $j\le d$, are pairwise disjoint subsets of $\frac32(B_\infty^d\cap V)$. Comparing
volumes in the $k$-dimensional space $V$ gives $d\,(\frac14)^k\le(\frac32)^k$, that is, $d\le6^k\le6^m$, a
contradiction. Hence there is $j$ with $\dist_\infty(e_j,V)\ge\frac14$, and by Lemma~\ref{lem:duality} there is
$\xi\in V^\perp$ with $\norm\xi_1\le1$ and $\xi_j=\ip{e_j}\xi\ge\frac14$.
\end{proof}

\begin{remark}\label{rem:ex-det}
Proposition~\ref{prop:ex-det} is known in a different language. The proof shows that a $d\times d$ matrix whose
diagonal entries differ from one by less than $\frac14$ and whose other entries are smaller than $\frac14$ in
absolute value has rank at least $\log_6d$. Lower bounds of this kind for the rank of perturbed identity matrices, also for small
perturbations, are due to Alon~\cite{Al09}, and the problem of recovering all coordinates of $\xi$ up to an error
$\eps\norm\xi_1$ from linear measurements is studied in the literature on data streams, see \cite{NNW14} for
deterministic algorithms. The exponential dimension in Proposition~\ref{prop:ex-det} is necessary: a number of
measurements of order $\eps^{-2}\log d$ suffices to recover all coordinates of $\xi$ up to an error $\eps\norm\xi_1$,
see \cite{NNW14}.
\end{remark}

\subsection{Randomized algorithms}
The algorithm behind Theorem~\ref{thm:lower} has three ingredients: a random hashing of the coordinates of
$\xi$, which removes the dependence on the dimension; a search for the large coordinates by repeated random
hashing; and a
Monte Carlo estimate of $\norm\xi_1$ with the help of Cauchy variables as in \cite{In06}. We start with two lemmas
on the first two ingredients. By hashing we mean that the coordinates are distributed randomly into groups, and that the measurements are the
sums of the coordinates over these groups.

In Proposition~\ref{prop:ex-det} the dimension $d$ is exponential in the number $m$ of measurements of the
deterministic algorithm. The search for the large coordinates in Lemma~\ref{lem:sketch} below needs a number of
measurements that grows with the logarithm of the dimension, and for $d>6^m$ this is at least of order $m$. But we want a
bound for all $m$. We therefore reduce the dimension first, by
hashing the coordinates of $\xi$ into $M$ groups. This
does not preserve $\norm\xi_\infty$, but it almost preserves the quantity \eqref{eq:S-example} that we are
interested in.

\begin{lemma}[Random hashing]\label{lem:hash}
Let $\xi\in\R^d$, $\tau\in(0,1]$ and $M\in\N$ with $M\tau\ge4$. Let $h(1),\dots,h(d)$ be independent random
variables that are uniformly distributed on $\{1,\dots,M\}$, and let $z\in\R^M$ be given by
$z_b:=\sum_{j\colon h(j)=b}\xi_j$. Then $\norm z_1\le\norm\xi_1$ and
\[
 \E\,\Bigl|\max\bigl\{\norm z_\infty,\tau\norm\xi_1\bigr\}-\max\bigl\{\norm\xi_\infty,\tau\norm\xi_1\bigr\}\Bigr|
 \le\frac{80}{\tau^2M}\,\norm\xi_1 .
\]
\end{lemma}

\begin{proof}
The first statement is clear, and the second is trivial for $\xi=0$. Both sides of the inequality are homogeneous
in $\xi$, so we may assume
$\norm\xi_1=1$. Let $\theta:=\tau/4$. We call a coordinate $j$ \emph{heavy} if $\abs{\xi_j}>\theta$, and
\emph{light} otherwise. Let $H$ be the set of heavy coordinates, so that $\abs H<1/\theta$. For $b\le M$ let
\[
 \Lambda_b:=\sum_{j\notin H\colon h(j)=b}\abs{\xi_j}
\]
be the light mass of the group $b$, and let $R:=\sum_{j\in H}\Lambda_{h(j)}$. We consider the events
\[
 E_1:=\{h\text{ is injective on }H\}\qquad\text{and}\qquad E_2:=\{\Lambda_b\le2\theta\text{ for all }b\le M\}.
\]

\smallskip
\emph{On $E_1\cap E_2$ the difference of the two maxima is at most $R$.} Indeed, a group $b$ without a heavy
coordinate satisfies $\abs{z_b}\le\Lambda_b\le2\theta\le\tau$. A group $b$ with a heavy coordinate $j$ contains no
other heavy coordinate, so that $\abs{z_b-\xi_j}\le\Lambda_b\le R$. This shows
$\norm z_\infty\le\max\{\tau,\norm\xi_\infty+R\}$. On the other hand, if $\norm\xi_\infty>\theta$, then a coordinate
with $\abs{\xi_j}=\norm\xi_\infty$ is heavy and hence $\norm z_\infty\ge\norm\xi_\infty-R$; and if
$\norm\xi_\infty\le\theta$, then $\max\{\norm\xi_\infty,\tau\}=\tau$.

Since both maxima are in $[\tau,1]$, their difference is at most one in any case, and it remains to show that
$\E R+\P(E_1^c)+\P(E_2^c)\le5/(\theta^2M)$.

\smallskip
\emph{The mean of $R$.} The variables $\Lambda_b$ depend only on $h(i)$ for light $i$, and are therefore independent
of $h(j)$ for heavy $j$. Since $\E\Lambda_b\le1/M$ for every $b$, we get $\E R\le\abs H/M\le1/(\theta M)$.

\smallskip
\emph{The event $E_1$.} Two given heavy coordinates are in the same group with probability $1/M$, hence
$\P(E_1^c)\le\binom{\abs H}2\frac1M\le1/(2\theta^2M)$.

\smallskip
\emph{The event $E_2$.} Fix $b$ and write $\Lambda_b=\sum_{j\notin H}c_j\beta_j$ with $c_j:=\abs{\xi_j}\le\theta$ and
independent Bernoulli variables $\beta_j$ with $\P(\beta_j=1)=1/M$. Let $\lambda:=\theta^{-1}\log(M\theta)\ge0$.
By convexity, $e^{\lambda c}-1\le\frac c\theta(e^{\lambda\theta}-1)\le cM$ for $0\le c\le\theta$. Hence
\[
 \E\,e^{\lambda\Lambda_b}=\prod_{j\notin H}\Bigl(1+\frac{e^{\lambda c_j}-1}M\Bigr)\le\prod_{j\notin H}(1+c_j)\le e,
\]
and Markov's inequality gives $\P(\Lambda_b>2\theta)\le e\cdot e^{-2\theta\lambda}=e\,(M\theta)^{-2}$. By the union
bound, $\P(E_2^c)\le e/(M\theta^2)$.

Adding the three bounds and using $\theta\le1$, we obtain the bound $(1+\frac12+e)/(\theta^2M)\le80/(\tau^2M)$.
\end{proof}

The second lemma is a standard tool for finding the large coordinates of a vector. The procedure is the count-min
sketch of Cormode and Muthukrishnan with a median in place of the minimum, see \cite[Theorem~2]{CM05}, and medians
of sums with random signs over random groups were used before in \cite{CCF04}; we include the short proof. Here and in the following, the
\emph{median} of real numbers $a_1,\dots,a_t$ is the $\lceil t/2\rceil$-th smallest of them. If $\mu$ is the
median, then at least $t/2$ of the numbers are at most $\mu$, and at least $t/2$ of them are at least $\mu$.

\begin{lemma}[Search by repeated hashing]\label{lem:sketch}
Let $M,t,w\in\N$ and $z\in\R^M$. Let $\chi_1,\dots,\chi_t$ be independent random mappings from $\{1,\dots,M\}$ to
$\{1,\dots,w\}$ whose values $\chi_i(b)$, $b\le M$, are independent and uniformly distributed on $\{1,\dots,w\}$. Let
\[
 s_{i,v}:=\sum_{b\colon\chi_i(b)=v}z_b\qquad\text{for }i\le t\text{ and }v\le w,
\]
and for $b\le M$ let $\hat z_b$ be the median of the numbers $s_{1,\chi_1(b)},\dots,s_{t,\chi_t(b)}$. Then
\[
 \P\Bigl(\norm{\hat z-z}_\infty>\frac4w\,\norm z_1\Bigr)\le M\,e^{-t/8}.
\]
\end{lemma}

So the $tw$ numbers $s_{i,v}$, each of which is a linear measurement of $z$, determine all coordinates of $z$ up to
an error of $4\norm z_1/w$, with high probability if $t$ is a sufficiently large multiple of $\log M$.\looseness=-1

\begin{proof}
We may assume $z\ne0$. Fix $b\le M$ and let $r:=4\norm z_1/w$. For $i\le t$ we have
\[
 s_{i,\chi_i(b)}-z_b=\sum_{b'\ne b}z_{b'}\,\mathbf 1\{\chi_i(b')=\chi_i(b)\},
\]
and $\P(\chi_i(b')=\chi_i(b))=1/w$ for $b'\ne b$. Hence $\E\abs{s_{i,\chi_i(b)}-z_b}\le\norm z_1/w$, and by Markov's
inequality the event $E_i:=\{\abs{s_{i,\chi_i(b)}-z_b}>r\}$ has probability at most $\frac14$; note that $r>0$. The events
$E_1,\dots,E_t$ are independent, since $E_i$ depends only on $\chi_i$. If $\hat z_b>z_b+r$, then at least $t/2$ of the numbers $s_{i,\chi_i(b)}$, $i\le t$,
are larger than $z_b+r$, and if $\hat z_b<z_b-r$, then at least $t/2$ of them are smaller than $z_b-r$. So
$\abs{\hat z_b-z_b}>r$ implies that at least $t/2$ of the events $E_i$ occur. Then the number of these events
exceeds its mean by at least $t/4$, and by Hoeffding's inequality for the indicators of the events $E_i$ this has
probability at most $\exp(-2(t/4)^2/t)=e^{-t/8}$. The union bound over $b\le M$ gives the claim.
\end{proof}

\begin{proposition}\label{prop:ex-ran}
There are absolute constants $C$ and $k_0$ such that for all $k\ge k_0$ and all $d\in\N$, with
\[
 \tau_k:=36\,\Bigl\lfloor\frac k{\lceil48\log(k+1)\rceil}\Bigr\rfloor^{-1}\le1,
\]
we have
\[
 e_{3k}^{\ran}(p_{\tau_k},B_1^d,\Lambda^{\rm all})\le C\,\frac{\log(k+1)}{k^{3/2}} .
\]
\end{proposition}

We describe the algorithm. The input is $\xi\in\R^d$ with $\norm\xi_1\le1$. The algorithm uses the parameters
\[
 t_k:=\lceil48\log(k+1)\rceil,\qquad w_k:=\Bigl\lfloor\frac k{t_k}\Bigr\rfloor\qquad\text{and}\qquad\eps_k:=\frac4{w_k},
\]
so that $\tau_k=9\,\eps_k$.
The algorithm has four stages, which use independent random numbers.
\begin{enumerate}
\item[\textbf{0.}] \emph{Reduction of the dimension.} Let $M:=(k+1)^4$, let $h$ be as in Lemma~\ref{lem:hash}, and
let $z\in\R^M$ with $z_b=\sum_{j\colon h(j)=b}\xi_j$. Every linear measurement of $z$ is a linear measurement of
$\xi$.
\item[\textbf{1.}] \emph{The $\ell_1$-norm.} Let $c_1,\dots,c_k\in\R^d$ be independent random vectors whose entries are
independent standard Cauchy variables. We measure $y_i=\ip{c_i}\xi$ and put
\[
 \hat\lambda:=\operatorname{median}\{\abs{y_1},\dots,\abs{y_k}\},\qquad\tilde\lambda:=\min\{\hat\lambda,1\}.
\]
\item[\textbf{2.}] \emph{Large coordinates.} Let $\chi_1,\dots,\chi_{t_k}$ and $\hat z\in\R^M$ be as in
Lemma~\ref{lem:sketch} with $t=t_k$ and $w=w_k$. That is, we measure the $t_kw_k\le k$ numbers
$s_{i,v}=\sum_{b\colon\chi_i(b)=v}z_b$ and compute the medians $\hat z_b$ for all $b\le M$. Let
\[
 T:=\bigl\{b\le M\colon\abs{\hat z_b}\ge4\,\eps_k\,\hat\lambda\bigr\}.
\]
If $1\le\abs T\le k$, we measure the coordinates $z_b$, $b\in T$, and put $a_T:=\max_{b\in T}\abs{z_b}$.
Otherwise we put $a_T:=0$.
\item[\textbf{3.}] \emph{Output.} We return $\hat S:=\max\{a_T,\ \tau_k\tilde\lambda\}$.
\end{enumerate}
Note that the $\ell_1$-norm in stage 1 is the norm of $\xi$ and not of $z$, while stage 2 works with $z$. The
threshold in the definition of $T$ is proportional to the estimate $\hat\lambda$, since we need an error that is
small relative to $\norm\xi_1$. All measurements are chosen in advance, except for the coordinates $z_b$, $b\in T$,
at the end of stage 2: the algorithm is adaptive only in its last step. It is measurable: the measurement vectors that are used before the last step are measurable functions of the
random numbers and do not depend on the observed
values, and the set $T$ and the output are obtained from the random numbers and the observed values by medians,
maxima and comparisons.

\begin{proof}[Proof of Proposition~\ref{prop:ex-ran}]
We choose $k_0$ such that $w_k\ge36$, and hence $\tau_k\le1$, for all $k\ge k_0$; one may take
$k_0=20\,000$. For $k\ge k_0$ we have $k\ge2\,t_k$ and $t_k\le49\log(k+1)$, hence $k/(2t_k)\le w_k\le k/t_k$ and
\begin{equation}\label{eq:tau-order}
 \frac9k\le\frac{1728\,\log(k+1)}k\le\tau_k=\frac{36}{w_k}\le\frac{72\,t_k}k\le\frac{3528\,\log(k+1)}k .
\end{equation}
We may assume $\xi\ne0$; for $\xi=0$ all observed values vanish, the set $T=\{1,\dots,M\}$ has more than $k$
elements, and the output is $0=p_{\tau_k}(0)$. Let
\[
 S:=p_{\tau_k}(\xi)=\max\{\norm\xi_\infty,\tau_k\norm\xi_1\}\qquad\text{and}\qquad S_z:=\max\{\norm z_\infty,\tau_k\norm\xi_1\}.
\]
Since $\norm z_\infty\le\norm z_1\le\norm\xi_1\le1$, the numbers $S$, $S_z$ and $\hat S$ are in $[0,1]$.

\smallskip
\emph{The number of measurements.} Stage 1 uses $k$ measurements. Stage 2 uses $t_kw_k\le k$ measurements for the
numbers $s_{i,v}$, and then at most $k$ coordinates are measured. Altogether, the algorithm uses at most $3k$ linear
measurements, and we may add superfluous measurements to obtain exactly $3k$.

\smallskip
\emph{Stage 0.} We have $M\tau_k\ge4$ and $\tau_k\ge9/k$ by \eqref{eq:tau-order}, so Lemma~\ref{lem:hash} gives
\begin{equation}\label{eq:hash-error}
 \E\abs{S_z-S}\le\frac{80}{\tau_k^2M}\le\frac{80\,k^2}{81\,(k+1)^4}\le\frac1{k^2}.
\end{equation}

\smallskip
\emph{Stage 1.} A linear combination $\sum_j\xi_jc_{ij}$ of independent standard Cauchy variables has the same
distribution as $\norm\xi_1$ times a standard Cauchy variable. Hence the variables $\abs{y_i}/\norm\xi_1$ are
independent with distribution function $x\mapsto\frac2\pi\arctan x$, $x\ge0$, whose median is $1$. For $u\in[0,1]$
we have $\frac2\pi\arctan(1+u)\ge\frac12+\frac{2u}{5\pi}$ and $\frac2\pi\arctan(1-u)\le\frac12-\frac u\pi$, because
the derivative of $\arctan$ is at least $\frac15$ on $[1,2]$ and at least $\frac12$ on $[0,1]$. If the median of
the $\abs{y_i}/\norm\xi_1$ exceeds $1+u$, then at least half of them exceed $1+u$, and similarly for $1-u$. By
Hoeffding's inequality,
\begin{equation}\label{eq:median}
 \P\Bigl(\bigl|\hat\lambda/\norm\xi_1-1\bigr|>u\Bigr)\le2\exp\bigl(-ku^2/31\bigr)\qquad\text{for }u\in[0,1].
\end{equation}
In particular, the event $G_1:=\{\frac12\norm\xi_1\le\hat\lambda\le2\norm\xi_1\}$ satisfies
$\P(G_1^c)\le2e^{-k/124}$. Moreover, since $\norm\xi_1\le1$, we have
$\abs{\tilde\lambda-\norm\xi_1}\le\min\{\abs{\hat\lambda/\norm\xi_1-1},1\}$ and therefore, by \eqref{eq:median},
\begin{equation}\label{eq:mass}
 \E\bigl|\tilde\lambda-\norm\xi_1\bigr|\le\int_0^12e^{-ku^2/31}\,du\le\sqrt{\frac{31\pi}k}\le\frac{10}{\sqrt k}.
\end{equation}
This paragraph does not use $k\ge k_0$: the bounds \eqref{eq:median} and \eqref{eq:mass} hold for every $k\in\N$.

\smallskip
\emph{Stage 2.} We condition on $h$ and on $c_1,\dots,c_k$. Then $z$ and $\hat\lambda$ are fixed, and the mappings
$\chi_1,\dots,\chi_{t_k}$ are independent of them. Since $\norm z_1\le\norm\xi_1$, $M=(k+1)^4$ and
$t_k\ge48\log(k+1)$, Lemma~\ref{lem:sketch} shows that the event
$G_2:=\{\norm{\hat z-z}_\infty\le\eps_k\norm\xi_1\}$ satisfies
\[
 \P(G_2^c\mid h,c_1,\dots,c_k)\le M\,e^{-t_k/8}\le(k+1)^{4}\,(k+1)^{-6}=(k+1)^{-2}.
\]
On $G_1\cap G_2$ the following holds.
\begin{itemize}
\item If $b\in T$, then $\abs{\hat z_b}\ge4\,\eps_k\hat\lambda\ge2\,\eps_k\norm\xi_1$ and hence
$\abs{z_b}\ge\eps_k\norm\xi_1$. Since $\norm z_1\le\norm\xi_1$, this implies
$\abs T\le1/\eps_k=w_k/4\le k$. So the coordinates in $T$ are measured if $T\ne\emptyset$, and in any case $a_T\le\norm z_\infty$.
\item If $\abs{z_b}\ge\tau_k\norm\xi_1=9\,\eps_k\norm\xi_1$, then
$\abs{\hat z_b}\ge8\,\eps_k\norm\xi_1\ge4\,\eps_k\hat\lambda$ and hence $b\in T$.
\end{itemize}
Consequently, on $G_1\cap G_2$ we have $a_T\le\norm z_\infty$, with equality if $\norm z_\infty\ge\tau_k\norm\xi_1$,
since then a coordinate with $\abs{z_b}=\norm z_\infty$ belongs to $T$.

\smallskip
\emph{The error.} Let $a:=\norm z_\infty$, $\beta:=\tau_k\norm\xi_1$ and $\hat\beta:=\tau_k\tilde\lambda$, so that
$S_z=\max\{a,\beta\}$. On $G_1\cap G_2$ we argue as follows. If $a\ge\beta$, then $a_T=a$ and
$\abs{\hat S-S_z}=\max\{a,\hat\beta\}-a\le\abs{\hat\beta-\beta}$. If $a<\beta$, then $a_T<\beta$ and
$\abs{\hat S-S_z}=\abs{\max\{a_T,\hat\beta\}-\beta}\le\abs{\hat\beta-\beta}$. So in both cases
$\abs{\hat S-S_z}\le\tau_k\abs{\tilde\lambda-\norm\xi_1}$, while $\abs{\hat S-S_z}\le1$ always. The complement of
$G_1\cap G_2$ has probability at most $\P(G_1^c)+\P(G_2^c)\le2e^{-k/124}+(k+1)^{-2}$. With
\eqref{eq:hash-error}, \eqref{eq:mass} and \eqref{eq:tau-order} we obtain
\[
 \E\abs{\hat S-S}\le\E\abs{S_z-S}+\E\abs{\hat S-S_z}\le\frac1{k^2}+\frac{10\,\tau_k}{\sqrt k}+2e^{-k/124}+\frac1{(k+1)^{2}}
 \le C\,\frac{\log(k+1)}{k^{3/2}}. \qedhere
\]
\end{proof}

\begin{proof}[Proof of Theorem~\ref{thm:lower}]
Let $m,n\in\N$, and let $k_0$ be the number of Proposition~\ref{prop:ex-ran}. By \eqref{eq:translate}, the bounds of
\eqref{eq:ex-det-tau} and Proposition~\ref{prop:ex-ran} are bounds for $e_m^{\dt}(S,F_{d,\tau})$ and
$e_{3k}^{\ran}(S,F_{d,\tau_k})$. Note that
$e_n^{\ran}(S,F)\le e_N^{\ran}(S,F)$ for $N\le n$, since an algorithm may ignore function values. If $n<3k_0$, we take
any $\tau\in(0,1]$ and $d>6^m$ and use $e_m^{\dt}(S,F_{d,\tau})\ge\frac18$ and
$e_n^{\ran}(S,F_{d,\tau})\le1\le(3k_0)^{3/2}\,n^{-3/2}\log(n+1)/\log2$. Otherwise
let $k:=\lfloor n/3\rfloor$, so that $k\ge k_0$ and $k>\frac n3-1\ge\frac n9$, and let
$\tau:=\tau_k$ and
$d>6^m$. By \eqref{eq:ex-det-tau} and Proposition~\ref{prop:ex-ran}, whose constant $C$ we call $C_1$ here,
\[
 e_m^{\dt}(S,F_{d,\tau})\ge\frac18\qquad\text{and}\qquad
 e_n^{\ran}(S,F_{d,\tau})\le e_{3k}^{\ran}(S,F_{d,\tau})\le C_1\,\frac{\log(k+1)}{k^{3/2}}\le27\,C_1\,\frac{\log(n+1)}{n^{3/2}}. \qedhere
\]
\end{proof}

\begin{remark}\label{rem:example}
(1) The number $m$ of measurements of the deterministic algorithm enters only through the condition $d>6^m$,
and the bound for randomized algorithms does not depend on $d$. This is the purpose of stage~0. Without it,
Lemma~\ref{lem:sketch} would have to be applied in dimension $d$, with a number $t$ of mappings $\chi_i$ of order
$\log d$, which is at least of order $m$. One can also search for the large coordinates with the adaptive algorithms of
Kunsch and Wnuk \cite{KW}, whose number of measurements depends on the dimension only through $\log\log d$; adaptive
algorithms of this kind go back to Indyk, Price and Woodruff~\cite{IPW11}.

(2) \emph{Bernstein widths.} For $j<d$ we have
\[
 b_j(F_{d,\tau})=\max\Bigl\{\tau,\frac1{j+1}\Bigr\}.
\]
For the classes $F=F_{d,\tau_k}$ of Theorem~\ref{thm:lower} this means that $b_j(F)=\tau_k$ is of order $\log(n)/n$
for all $k\le j<d$, while $\rho_m(F)\ge\frac14$ by Proposition~\ref{prop:ex-det}. So the factor $n$ in
\eqref{eq:det-bern} is attained up to a logarithm if $m\ge8n$. Moreover,
Conjecture~\ref{conj:main} would give $e_n^{\ran}(S,F_{d,\tau_k})\ge c\,\tau_kn^{-1/2}$ for $d>Cn$, which matches
Proposition~\ref{prop:ex-ran} up to constants, while Theorem~\ref{thm:bernstein} only gives a lower bound of
order $\tau_k/n$.

To prove the formula, recall that, by Lemma~\ref{lem:dict}, $b_j(F_{d,\tau})$ is the supremum of all $\lambda\ge0$ for which there is a subspace
$U\subset\R^d$ of dimension $j+1$ such that $\norm\xi_1\le1$ for all $\xi\in U$ with $p_\tau(\xi)\le\lambda$.
The lower bound $\tau$ holds since $p_\tau(\xi)\le\tau$ implies $\norm\xi_1\le1$. The lower
bound $1/(j+1)$ is obtained for the space of all $\xi$ that are supported in the first $j+1$ coordinates, since
then $\norm\xi_1\le(j+1)\norm\xi_\infty\le(j+1)\,p_\tau(\xi)$. For the upper
bound let $U$ and $\lambda>\tau$ be as above. If $\xi\in U$ and $\norm\xi_\infty=\lambda$, then
$s:=p_\tau(\xi)\ge\lambda$ and $p_\tau(\frac\lambda s\xi)=\lambda$, hence $\norm\xi_1\le s/\lambda$. This
excludes $s=\tau\norm\xi_1$, hence $s=\lambda$ and $\norm\xi_1\le1$. So
$\lambda(B_\infty^d\cap U)\subset B_1^d$, and $\lambda\le1/(j+1)$ as in Remark~\ref{rem:sharp}.

(3) The logarithm in Theorem~\ref{thm:lower} comes from the number $t_k$ of mappings $\chi_i$ in stage~2, which
determines the failure probability of the median estimator in the search for the large coordinates.
\end{remark}

\begin{remark}[The max-norm and the $1$-norm]\label{rem:norms}
Parts of the algorithm give bounds for the estimation of the norms $\norm\cdot_\infty$ and $\norm\cdot_1$ on $B_1^d$.

(1) \emph{The max-norm.} For $k\ge k_0$ consider stage~0 and
the first part of stage~2 of the algorithm above, with the output $\min\{\norm{\hat z}_\infty,1\}$. This is a
randomized algorithm with at most $k$ linear measurements, which is even nonadaptive, and $\tau_k$ is now only a parameter of the
analysis. With $S$ and $S_z$ as in the proof of Proposition~\ref{prop:ex-ran}, we have
$\abs{\norm z_\infty-\norm\xi_\infty}\le\abs{S_z-S}+\tau_k$ and $\E\abs{S_z-S}\le k^{-2}$, see \eqref{eq:hash-error},
and on the event $G_2$ we have
$\abs{\min\{\norm{\hat z}_\infty,1\}-\norm z_\infty}\le\eps_k$. So the expected error is at
most $k^{-2}+\tau_k+\eps_k+(k+1)^{-2}$. Taking $k=n$ for $n\ge k_0$, and using the trivial bound $1$ for $n<k_0$, we
obtain
\[
 e_n^{\ran}(\norm\cdot_\infty,B_1^d,\Lambda^{\rm all})\le C\,\frac{\log(n+1)}n\qquad\text{for all }n,d\in\N .
\]
By Proposition~\ref{prop:ex-det}, this is an example where nonadaptive randomized algorithms are better than all
deterministic algorithms with $m$ measurements by a factor of order $n/\log n$, if $d>6^m$. By Lemma~\ref{lem:dict}, the corresponding problem for function values is the computation of
the supremum for the class $F_d$ of all $f_\xi=\ip\cdot\xi$ with $\norm\xi_1\le1$ on the domain $B_1^d$, where
$S(f_\xi)=\norm\xi_\infty$.

For the linear problem of approximating $\xi$ in the norm of
$\ell_\infty^d$ on the same set, which is the example of \cite{KNW,KW}, the hashing of Lemma~\ref{lem:hash} does not
help, since it loses the positions of the large coordinates. There, the error of nonadaptive randomized algorithms
necessarily depends on the dimension $d$ \cite{KNW}, and adaptive
randomized algorithms achieve an error of order $(\log n+\log\log d)/n$ \cite[Theorem~3.3]{KW}; it is not known
whether the dependence on the dimension $d$ can be removed for adaptive algorithms, see
the discussion after Open Problem~\ref{prob:point-value}.

(2) \emph{The $1$-norm.} Stage~1 of the algorithm alone, with the output $\tilde\lambda$, is a nonadaptive
randomized algorithm, and \eqref{eq:mass} shows that
\[
 e_n^{\ran}(\norm\cdot_1,B_1^d,\Lambda^{\rm all})\le\frac{10}{\sqrt n}\qquad\text{for all }n,d\in\N .
\]
By Lemma~\ref{lem:dict}, the corresponding problem for function values is the computation of the supremum for the
class $F_{d,1}$ of all $f_\xi=\ip\cdot\xi$ with $\norm\xi_1\le1$ on the cube $[-1,1]^d$, which is the unit ball of
the dual norm; here $S(f_\xi)=\norm\xi_1$. This class is the
unit ball of a $d$-dimensional subspace of $B([-1,1]^d)$, so that $b_j(F_{d,1})=1$ for all $j<d$. This shows that the
factor $n^{-1/2}$ in Conjecture~\ref{conj:main} would be optimal.
\end{remark}

\section{Finding a maximizer}\label{sec:opt}

Instead of the number $S(f)$ one may want to find a point of $D$ at which the supremum is almost attained. We
denote this problem by $\OPT$ and collect what the results of this paper imply for it. Since $F$ is symmetric,
everything applies to the problem of finding a minimizer as well.

A deterministic algorithm for $\OPT$ using $n$ function values is a mapping $A_n\colon F\to D$ of the form
$A_n(f)=\varphi(y_1,\dots,y_n)$, where $y_i$ and $\psi_i$ are as in Section~\ref{sec:prelim} and
$\varphi\colon\R^n\to D$. Its error is the largest residual,
\[
 e(A_n,\OPT,F):=\sup_{f\in F}\bigl(S(f)-f(A_n(f))\bigr),
\]
and the numbers $e_n^{\dt}(\OPT,F)$ and $e_n^{\ran}(\OPT,F)$ are defined as in
Section~\ref{sec:prelim}; the error of a randomized algorithm is $\sup_{f\in F}\E_\omega\bigl(S(f)-f(A_n^\omega(f))\bigr)$.
For measurable randomized algorithms we require in addition that the mappings
$(\omega,y)\mapsto\bigl(g_l(\varphi^\omega(y))\bigr)_{l\le r}$ are measurable, for the subspaces and bases considered
there.

\begin{proposition}\label{prop:opt}
Let $F\subset B(D)$ be nonempty, convex and symmetric, and $n\in\N$. Then
\[
 \tfrac12\,\rho_n(F)\;\le\;e_n^{\dt}(\OPT,F)\;\le\;2\,e_n^{\dt}(\APP,F)\;=\;2\,\rho_n(F)
\]
and
\[
 e_{n+1}^{\dt}(S,F)\le e_n^{\dt}(\OPT,F),\qquad e_{n+1}^{\ran}(S,F)\le e_n^{\ran}(\OPT,F),\qquad
 e_n^{\ran}(\OPT,F)\le2\,\rho_n(F).
\]
\end{proposition}

So, up to one function value and a factor two, computing the supremum is at most as hard as finding a maximizer,
which is at most as hard as uniform approximation; for deterministic algorithms the three problems are
equivalent. The deterministic part is essentially known: Novak~\cite[Proposition~1.3.2]{No88} proves
$\rho_{n+1}(F)\le e_n^{\dt}(\OPT,F)\le2\,\rho_n(F)$, and the first inequality of the second line is contained in
\cite[Theorem~5.4]{Wa84}. We include the short proofs. The last inequality
shows that $e_n^{\ran}(\OPT,F)\le4\,e_n^{\dt}(\OPT,F)$, see Remark~\ref{rem:ran-det}.

\begin{proof}
\emph{First inequality.} Let $A_n$ be a deterministic algorithm for $\OPT$, let $x_1,\dots,x_n$ be the points that it
uses for the input $f=0$, and let $h\in F$ vanish at these points. As in the proof of Proposition~\ref{prop:det}, the
algorithm returns the same point $\hat x$ for $h$ and for $-h$. The two residuals add up to
$S(h)-h(\hat x)+S(-h)+h(\hat x)=S(h)+S(-h)\ge\norm h_\infty$, since $S(h)$ and $S(-h)$ are nonnegative (as
$h(x_1)=0$) and their maximum is $\norm h_\infty$. Hence the error is at least $\frac12\norm h_\infty$.

\smallskip
\emph{Second inequality.} Let $A_n$ be a deterministic algorithm for uniform approximation with error $\eps<\infty$, and let
$\eta>0$. For an input $f$ let $g:=A_n(f)$, and return a point $\hat x\in D$ with $g(\hat x)\ge S(g)-\eta$; this point
depends only on the observed function values. Since $\norm{f-g}_\infty\le\eps$, we get
$S(f)-f(\hat x)\le S(g)+\eps-g(\hat x)+\eps\le2\eps+\eta$. The equality in the first line was noted in
Remark~\ref{rem:app}.

\smallskip
\emph{Second line.} If an algorithm for $\OPT$ returns the point $\hat x$, then one further function value gives the
number $f(\hat x)$, and $\abs{S(f)-f(\hat x)}$ is the residual. This turns deterministic and measurable randomized
algorithms for $\OPT$ with $n$ function values into algorithms for $S$ with $n+1$ function values and the same error.

\smallskip
\emph{Last inequality.} We may assume that $\rho_n(F)$ is finite. Let $x_1,\dots,x_n$, $\rho$, $g_y$ and $L$ be as in
the proof of Proposition~\ref{prop:meas}, and let $\eta>0$. We choose a partition of $\R^n$ into countably many
Borel sets $Q_1,Q_2,\dots$ such that $L\,\norm{y-y'}_1\le\eta$ whenever $y$ and $y'$ belong to the same set. For every
$m$ we choose $y_m\in Q_m$ and $z_m\in D$ with $g_{y_m}(z_m)\ge S(g_{y_m})-\eta$, and we put $\varphi(y):=z_m$ for
$y\in Q_m$. Then, for $y\in Q_m$,
\[
 g_y(\varphi(y))\ge g_{y_m}(z_m)-\eta\ge S(g_{y_m})-2\eta\ge S(g_y)-3\eta .
\]
For $f\in F$ and $y=(f(x_1),\dots,f(x_n))$ we have $\norm{f-g_y}_\infty\le\rho$ and therefore
$S(f)-f(\varphi(y))\le S(g_y)+\rho-g_y(\varphi(y))+\rho\le2\rho+3\eta$. So the nonadaptive algorithm
$f\mapsto\varphi(y)$ has error at most $2\rho+3\eta$. It is a measurable randomized algorithm, which does not
depend on $\omega$: the points are fixed, and for every function $g$ on $D$ the mapping $y\mapsto g(\varphi(y))$ is
constant on each of the sets $Q_m$ and hence measurable.
\end{proof}

Together with \eqref{eq:app-upper}, the two lines of the proposition give the upper bound of Theorem~\ref{thm:upper} for $\OPT$.

\begin{corollary}\label{cor:opt-upper}
Let $F\subset B(D)$ be nonempty, convex and symmetric. Then, for all $n\ge50$ and with $C=22\,400$,
\[
 e_{8n-1}^{\dt}(\OPT,F)\;\le\;2\,C\,n^2\,e_{n-1}^{\ran}(\OPT,F).
\]
\end{corollary}

In the same way, Conjecture~\ref{conj:main} would give an upper bound of order $n^{3/2}$ for $\OPT$. One might hope
for more: for deterministic algorithms $\OPT$ is equivalent to uniform approximation, and for uniform approximation
Bernstein widths are lower bounds for randomized algorithms without a loss in $n$, so that the gain of randomization
is at most of order $n$, see Remark~\ref{rem:app-upper}. The following example shows that $\OPT$ is different. As
for $S$, randomized algorithms can beat the Bernstein widths by a factor of order $n^{-1/2}$.

\begin{proposition}\label{prop:block}
Let $C\in\N$ and $C':=60\,C^{3/2}$. For all $n\ge4C$ there is a set $D$ and a nonempty, bounded, convex and symmetric class
$F\subset B(D)$ with
\[
 e_n^{\ran}(\OPT,F)\;\le\;C'\,n^{-1/2}\;b_{Cn-1}(F).
\]
\end{proposition}

\begin{proof}
Let $M,t,k\in\N$ and $d:=Mt$. We write $\xi\in\R^d$ as $\xi=(\xi_1,\dots,\xi_M)$ with blocks $\xi_b\in\R^t$. Let
$Q_b$ be the set of all $x\in\R^d$ with $x_b\in[-1,1]^t$ and $x_{b'}=0$ for $b'\ne b$, let $D:=\bigcup_{b\le M}Q_b$,
and let $F$ be the class of all $f_\xi=\ip\cdot\xi$ on $D$ with $\max_{b\le M}\norm{\xi_b}_1\le1$.
We have $S(f_\xi)=\norm{f_\xi}_\infty=\max_b\norm{\xi_b}_1$ for all $\xi\in\R^d$. So $F$ is the unit ball of a
$d$-dimensional subspace of $B(D)$, and $b_j(F)=1$ for $j<d$. We show that
\begin{equation}\label{eq:block}
 e_{Mk+t}^{\ran}(\OPT,F)\;\le\;\frac{10\,M}{\sqrt k}.
\end{equation}
For $x\in Q_b$ we have $f_\xi(x)=\ip{x_b}{\xi_b}$, so
every linear measurement of a single block $\xi_b$ is obtained from one function value, as in
Lemma~\ref{lem:dict}. The algorithm first applies stage~1 of the algorithm of Section~\ref{sec:example} to each
block $\xi_b$, with $k$ Cauchy measurements per block, and obtains numbers $\tilde\lambda_b\in[0,1]$ with
$\E\abs{\tilde\lambda_b-\norm{\xi_b}_1}\le10/\sqrt k$, see \eqref{eq:mass}; for $\xi_b=0$ we have
$\tilde\lambda_b=0$. Then it chooses the smallest $\hat b$ with
$\tilde\lambda_{\hat b}=\max_b\tilde\lambda_b$, measures the $t$ coordinates of $\xi_{\hat b}$, and returns the point
$\hat x\in Q_{\hat b}$ whose coordinates in the block $\hat b$ are the signs of the coordinates of $\xi_{\hat b}$.
Then $f_\xi(\hat x)=\norm{\xi_{\hat b}}_1$. Let $b^*$ be a block with $\norm{\xi_{b^*}}_1=S(f_\xi)$. If $\hat b=b^*$,
the residual is zero. Otherwise $\tilde\lambda_{b^*}\le\tilde\lambda_{\hat b}$ gives
\[
 S(f_\xi)-f_\xi(\hat x)\le\bigl(\norm{\xi_{b^*}}_1-\tilde\lambda_{b^*}\bigr)+\bigl(\tilde\lambda_{\hat b}-\norm{\xi_{\hat b}}_1\bigr)
 \le\sum_{b\le M}\bigl|\tilde\lambda_b-\norm{\xi_b}_1\bigr| ,
\]
and the expectation of the right-hand side is at most $10M/\sqrt k$. The algorithm uses $Mk+t$ function values, and
it is measurable for the same reasons as the algorithm of Section~\ref{sec:example}. This proves \eqref{eq:block}. Now let
$t:=\lceil n/2\rceil$, $M:=2C$ and $k:=\lfloor n/(4C)\rfloor$. Then $Mk+t\le n$, since $Mk$ is an integer and at most $\frac n2$. Moreover, $d\ge Cn$, so that $b_{Cn-1}(F)=1$, and
$k\ge n/(8C)$. Since an algorithm may ignore function values, \eqref{eq:block} gives
\[
 e_n^{\ran}(\OPT,F)\le e_{Mk+t}^{\ran}(\OPT,F)\le\frac{20\,C}{\sqrt k}\le20\sqrt8\,C^{3/2}\,n^{-1/2}\le C'\,n^{-1/2}\,b_{Cn-1}(F).
\]
Note that the algorithm uses both adaption and randomization.
\end{proof}

In the other direction, that is, for lower bounds on the gain of randomization for $\OPT$, we only have the following
consequence of a result of Kunsch and Wnuk~\cite{KW}.

\begin{proposition}\label{prop:opt-lower}
There is an absolute constant $c>0$ such that for all $m,n\in\N$ there is a set $D$ and a nonempty, bounded, convex
and symmetric class $F\subset B(D)$ with
\[
 e_m^{\dt}(\OPT,F)\;\ge\;c\,\frac{n}{\log(n+1)+\log(m+1)}\;e_n^{\ran}(\OPT,F).
\]
\end{proposition}

\begin{proof}
Let $F=F_d$ be the class of Remark~\ref{rem:norms}(1) on the domain $B_1^d$ with $d:=6^m+16$. Proposition~\ref{prop:ex-det} shows $\rho_m(F_d)\ge\frac14$, hence $e_m^{\dt}(\OPT,F_d)\ge\frac18$ by
Proposition~\ref{prop:opt}. By \cite[Theorem~3.3]{KW} there is a randomized algorithm that uses $n$ linear
measurements of $\xi$ and returns $\hat\xi\in\R^d$ with
$\E\norm{\xi-\hat\xi}_\infty\le C\,(\log(n+1)+\log\log d)/n$ whenever $\norm\xi_1\le1$. (The result is stated there up to constants and for large $d$, and its proof uses
$d\ge16$. For small $n$ the bound holds trivially if $C$ is large, since $\hat\xi=0$ has error at most one, and the
error is nondecreasing in $d$, since a vector can be padded with zeros.) This algorithm is described in \cite[Section~3.1]{KW}. It is measurable in the sense
above: it uses finitely many random signs and random hashings of the coordinates, each measurement vector is
chosen from a finite set, and the choice depends on the values observed before only through comparisons of medians
of their absolute values \cite[Section~2.2]{KW} and through rounded quotients of two of them
\cite[Section~2.3]{KW}, where ties are broken by a fixed rule; the
output $\hat\xi$ consists of observed coordinates and zeros.
Choose the smallest $j$ with $\abs{\hat\xi_j}=\norm{\hat\xi}_\infty$ and return the point $\hat x$
which is $e_j$ if $\hat\xi_j\ge0$ and $-e_j$ otherwise. Then
\[
 f_\xi(\hat x)\ge\abs{\hat\xi_j}-\abs{\xi_j-\hat\xi_j}\ge\norm\xi_\infty-2\,\norm{\xi-\hat\xi}_\infty ,
\]
so that the expected residual is at most $2C\,(\log(n+1)+\log\log d)/n$, and $\log\log d\le3\log(m+1)$.
\end{proof}

So randomization can help for $\OPT$ by a factor of order $n/\log n$ if $m$ is at most polynomial in
$n$. In contrast to Theorem~\ref{thm:lower}, the bound depends on $m$: the random hashing of
Lemma~\ref{lem:hash} removes the dimension from the problem of computing the largest coordinate up to a small
error, but not from the problem of locating it. For the classes of Theorem~\ref{thm:lower} we do not get more than
this either, since the algorithm of Section~\ref{sec:example} does not find a point at which the term
$\tau\norm\xi_1$ is almost attained. In fact, for randomized algorithms, unlike for deterministic algorithms, finding a
maximizer can be much harder than computing the supremum.

\begin{theorem}\label{thm:point-value}
\begin{enumerate}
\item[(a)] Let $F\subset B(D)$ be nonempty, convex and symmetric. Then, for all $n\ge50$ and with $C=22\,400$,
\[
 e_{8n-1}^{\ran}(\OPT,F)\;\le\;2\,C\,n^2\;e_n^{\ran}(S,F).
\]
\item[(b)] There is an absolute constant $c>0$ such that for all $m,n\in\N$ there is a set $D$ and a nonempty, bounded, convex
and symmetric class $F\subset B(D)$ with
\[
 e_m^{\ran}(\OPT,F)\;\ge\;c\,\frac{n}{\log(n+1)}\;e_n^{\ran}(S,F).
\]
More precisely, $e_m^{\ran}(\OPT,F)\ge\frac14$ and $e_n^{\ran}(S,F)\le C'\log(n+1)/n$ with an absolute constant $C'$.
\end{enumerate}
\end{theorem}

\begin{proof}
(a) By Proposition~\ref{prop:opt}, Remark~\ref{rem:app} and \eqref{eq:app-upper},
\[
 e_{8n-1}^{\ran}(\OPT,F)\le2\,\rho_{8n-1}(F)\le2\,C\,n^2\,e_n^{\ran}(S,F,\Lambda^{\rm all})
 \le2\,C\,n^2\,e_n^{\ran}(S,F).
\]

\smallskip
(b) Let $\tau:=1/(n+1)$ and $d>6^m$, so that $d\ge2\cdot3^m$, let
\[
 D:=\{\pm e_j\colon j\le d\}\cup\tau\,\{-1,0,1\}^d\subset\R^d ,
\]
and let $F$ be the class of all $f_\xi=\ip\cdot\xi$ on $D$ with $\norm\xi_1\le1$. Then
$S(f_\xi)=\norm{f_\xi}_\infty=\max\{\norm\xi_\infty,\tau\norm\xi_1\}$, which is the norm $p_\tau(\xi)$ from \eqref{eq:S-example}. The difference to Section~\ref{sec:example} is that $D$ is not
the unit ball of the dual norm, so that not all linear measurements are available: the available measurements are the coordinates of $\xi$ and the numbers $\ip w\xi$ with
$w\in\{-1,0,1\}^d$.

\smallskip
\emph{Lower bound.} Consider the inputs $f_{e_J}$, $J\le d$. On these inputs the function value at any given point
of $D$ takes at most three different values, since $f_{e_J}(\pm e_j)\in\{0,\pm1\}$ and
$f_{e_J}(\tau w)=\tau w_J\in\{-\tau,0,\tau\}$.
Hence a deterministic algorithm for $\OPT$ with $m$ function values observes at most $3^m$ different vectors of
values on these inputs, since each point is determined by the values observed before, and it returns at most $3^m$
different points. If the returned point is not $e_J$, then the
residual for $f_{e_J}$ is at least $1-\tau$. So the average residual with respect to the uniform distribution on
$\{f_{e_J}\colon J\le d\}$ is at least $(1-\tau)(1-3^m/d)\ge\frac14$, and the same holds for randomized algorithms
by the averaging argument \eqref{eq:bakhvalov}. This argument applies in the same way to the residual
$S(f)-f(A^\omega(f))$ in place of the error $\abs{T(f)-A^\omega(f)}$, since the residual is jointly measurable by the
additional requirement for $\OPT$.

\smallskip
\emph{Upper bound.} The algorithm of Remark~\ref{rem:norms}(1) only uses measurements of the form
$\sum_{j\in A}\xi_j$ with sets $A\subset\{1,\dots,d\}$, which are obtained from the function values at the points
$\tau\mathbf 1_A\in D$. Its analysis does not depend on the domain. So, if $n\ge k_0$, it computes $\norm\xi_\infty$
with $n$ function values and an expected error of at most $C_1\log(n+1)/n$, where $C_1$ denotes the constant $C$ of
Remark~\ref{rem:norms}(1). Since
$0\le S(f_\xi)-\norm\xi_\infty\le\tau=1/(n+1)$, this gives $e_n^{\ran}(S,F)\le(C_1+2)\log(n+1)/n$. If $n<k_0$, we use
that $e_n^{\ran}(S,F)\le1$, since $0\le S(f)\le1$ for all $f\in F$.
\end{proof}

On the classes of Theorem~\ref{thm:point-value}(b) randomization does not help for $\OPT$, up to a constant,
since $e_m^{\dt}(\OPT,F)\le2\rho_m(F)\le2$, while it helps for $S$ by a factor of order $n/\log n$, since
$e_m^{\dt}(S,F)\ge\frac18$; this follows from Propositions~\ref{prop:ex-det} and~\ref{prop:det}, since $d>6^m$ and the unit vectors
$\pm e_j$ belong to $D$, so that $\rho_m(F)\ge\frac14$. Let us also note that the bound of
Theorem~\ref{thm:point-value}(a) holds for $e_{8n-1}^{\dt}(\OPT,F)$ as well, with the same proof: finding a maximizer
is harder than computing the supremum by a factor of order $n^2$ at most, even if only the latter may use
randomization.

\section{Open problems}\label{sec:open}

We collect some questions that are left open.

\begin{problem}
Prove or disprove Conjecture~\ref{conj:main}, possibly with additional logarithmic factors. A proof would show that
$e_m^{\dt}(S,F)\le n^{3/2+o(1)}\,e_n^{\ran}(S,F)$ for all convex and symmetric classes and all $m\ge n^{1+o(1)}$,
which is optimal by Theorem~\ref{thm:lower}.
\end{problem}

\begin{problem}\label{prob:oversampling}
Is there a polynomial bound for $e_n^{\dt}(S,F)/e_n^{\ran}(S,F)$, that is, without oversampling? Our bounds require
that the deterministic algorithm uses at least $4n-1$ function values. We do not
know an example where the ratio $e_n^{\dt}(S,F)/e_n^{\ran}(S,F)$ is larger than $n^{3/2}$.
\end{problem}

\begin{problem}\label{prob:rates}
The classes in Theorem~\ref{thm:lower} depend on $n$. For a single class $F$, how large can the difference between
the polynomial rates of convergence of $e_n^{\ran}(S,F)$ and $e_n^{\dt}(S,F)$ be? Theorem~\ref{thm:upper} shows
that it is at most $2$. Theorem~\ref{thm:lower} does not give a lower bound.
We believe that a difference of $\frac54$ can be obtained by combining building blocks of the form of
Lemma~\ref{lem:dict}, with $B=B_1^d$ and the norm of $\ell_4^d$, for a sequence of dimensions $d$ into a single class.
\end{problem}

\begin{problem}\label{prob:opt-gain}
What is the maximal gain of randomization for the problem $\OPT$ of finding a maximizer? By
Corollary~\ref{cor:opt-upper}, the ratio $e^{\dt}_{8n-1}(\OPT,F)/e^{\ran}_{n-1}(\OPT,F)$ is at most of order $n^2$ for
every class, and Conjecture~\ref{conj:main} would improve this to $n^{3/2}$. By Proposition~\ref{prop:opt-lower},
the ratio $e^{\dt}_m(\OPT,F)/e^{\ran}_n(\OPT,F)$ can be of order $n/\log n$ if $m$ is at most polynomial in $n$. We do
not know a class for which it is larger than $n$.
\end{problem}

\begin{problem}\label{prob:point-value}
How much harder than computing the supremum can it be to find a maximizer, if randomized algorithms are
used for both problems? We always have $e^{\ran}_{n+1}(S,F)\le e^{\ran}_n(\OPT,F)$. By
Theorem~\ref{thm:point-value}, the ratio $e^{\ran}_{m}(\OPT,F)/e^{\ran}_n(S,F)$ can be of order $n/\log n$ for every
$m$, and it is at most of order $n^2$ for $m=8n-1$.
\end{problem}

A question related to Open Problem~\ref{prob:point-value} is
whether the largest coordinate of a vector $\xi\in\R^d$ with $\norm\xi_1\le1$ can be located, and not only computed,
by a randomized adaptive algorithm with a number of linear measurements that does not depend on $d$; the results of
\cite{KW} depend on $d$ through $\log\log d$. This is a question from adaptive sparse recovery, see \cite{IPW11}.
A number of measurements of order $\log\log d$ is known to be necessary for sparse recovery with an error measured
in the Euclidean norm \cite{PW13}. In the hard instances there, the $\ell_1$-norm of the vector exceeds its largest
coordinate by a factor of order $\sqrt d$, so that this lower bound does not apply to vectors with
$\norm\xi_1\le1$, and the optimal dependence on the dimension is stated as an open problem in \cite{KW2}. If this is not possible, in the sense that $e_n^{\ran}(\OPT,F_d)$ does
not tend to zero uniformly in $d$ for the classes $F_d$ of Remark~\ref{rem:norms}(1), then
$e^{\ran}_n(\OPT,F)/e^{\ran}_n(S,F)$ would be at least of order $n^{3/2}/\log n$ for the classes of Theorem~\ref{thm:lower}
with very large $d$.

\section*{Acknowledgement}
The AI models Claude Fable~5.1, Claude Fable~5.2 and Claude Opus~5.5 (Anthropic) have been used to assist the author in all aspects of this
work. The author retains full responsibility for the content of this paper.


\begin{thebibliography}{99}

\bibitem{Al09} N.~Alon, Perturbed identity matrices have high rank: proof and applications, Combin.\ Probab.\
Comput.\ 18 (2009), 3--15.

\bibitem{Ba59} N.~S.~Bakhvalov, On the approximate calculation of multiple integrals, Vestnik MGU, Ser.\ Math.\
Mech.\ Astron.\ Phys.\ Chem.\ 4 (1959), 3--18; English translation: J.~Complexity 31 (2015), 502--516.

\bibitem{Ba71} N.~S.~Bakhvalov, On the optimality of linear methods for operator approximation in convex classes of
functions, USSR Comput.\ Math.\ Math.\ Phys.\ 11 (1971), 244--249.

\bibitem{CCF04} M.~Charikar, K.~Chen, M.~Farach-Colton, Finding frequent items in data streams, Theoret.\ Comput.\
Sci.\ 312 (2004), 3--15.

\bibitem{CM05} G.~Cormode, S.~Muthukrishnan, An improved data stream summary: the count-min sketch and its
applications, J.~Algorithms 55 (2005), 58--75.

\bibitem{CW04} J.~Creutzig, P.~Wojtaszczyk, Linear vs.\ nonlinear algorithms for linear problems, J.~Complexity 20
(2004), 807--820.

\bibitem{Ga02} R.~J.~Gardner, The Brunn--Minkowski inequality, Bull.\ Amer.\ Math.\ Soc.\ 39 (2002), 355--405.

\bibitem{He92} S.~Heinrich, Lower bounds for the complexity of Monte Carlo function approximation,
J.~Complexity 8 (1992), 277--300.

\bibitem{In06} P.~Indyk, Stable distributions, pseudorandom generators, embeddings, and data stream computation,
J.~ACM 53 (2006), 307--323.

\bibitem{IPW11} P.~Indyk, E.~Price, D.~P.~Woodruff, On the power of adaptivity in sparse recovery, in: 52nd Annual
IEEE Symposium on Foundations of Computer Science (FOCS 2011), IEEE, 2011, 285--294.

\bibitem{KNU} D.~Krieg, E.~Novak, M.~Ullrich, On the power of adaption and randomization, Forum Math.\ Sigma 13
(2025), e152.

\bibitem{KNU2} D.~Krieg, E.~Novak, M.~Ullrich, On the power of adaption and randomization 2: non-symmetric convex
sets, preprint, 2026, arXiv:2610.33733.

\bibitem{KPUU} D.~Krieg, K.~Pozharska, M.~Ullrich, T.~Ullrich, Sampling projections in the uniform norm,
J.~Math.\ Anal.\ Appl.\ 553 (2026), no.~2, 129873, arXiv:2401.02220.

\bibitem{KU26} D.~Krieg, M.~Ullrich, Approximation of functions: optimal sampling and complexity, Acta Numer.\
35 (2026), 273--457, arXiv:2602.02066.

\bibitem{Ku16} R.~J.~Kunsch, Bernstein numbers and lower bounds for the Monte Carlo error, in: Monte Carlo and
Quasi-Monte Carlo Methods (MCQMC 2014), Springer Proc.\ Math.\ Stat.\ 163, Springer, 2016, 471--488.

\bibitem{Ku17} R.~J.~Kunsch, High-dimensional function approximation: breaking the curse with Monte Carlo
methods, Ph.D.\ dissertation, Friedrich Schiller University Jena, 2017, arXiv:1704.08213.

\bibitem{KNW} R.~J.~Kunsch, E.~Novak, M.~Wnuk, Randomized approximation of summable sequences -- adaptive and
non-adaptive, J.~Approx.\ Theory 304 (2024), 106056.

\bibitem{KW} R.~J.~Kunsch, M.~Wnuk, Uniform approximation of vectors using adaptive randomized information,
J.~Approx.\ Theory 313 (2026), 106216, arXiv:2408.01098.

\bibitem{KW2} R.~J.~Kunsch, M.~Wnuk, Adaptive and non-adaptive randomized approximation of high-dimensional
vectors, preprint, 2024, arXiv:2410.23067.

\bibitem{MH63} B.~S.~Mityagin, G.~M.~Henkin, Inequalities between $n$-diameters, in: Proc.\ Seminar on Functional
Analysis 7, Voronezh, 1963, 97--103 (in Russian).

\bibitem{NNW14} J.~Nelson, H.~L.~Nguyen, D.~P.~Woodruff, On deterministic sketching and streaming for sparse recovery
and norm estimation, Linear Algebra Appl.\ 441 (2014), 152--167.

\bibitem{NY} A.~S.~Nemirovsky, D.~B.~Yudin, Problem Complexity and Method Efficiency in Optimization, Wiley, 1983.

\bibitem{No88} E.~Novak, Deterministic and Stochastic Error Bounds in Numerical Analysis, Lecture Notes in Math.\
1349, Springer, 1988.

\bibitem{No92} E.~Novak, Optimal linear randomized methods for linear operators in Hilbert spaces,
J.~Complexity 8 (1992), 22--36.

\bibitem{No96} E.~Novak, On the power of adaption, J.~Complexity 12 (1996), 199--237.

\bibitem{PW13} E.~Price, D.~P.~Woodruff, Lower bounds for adaptive sparse recovery, in: Proc.\ 24th Annual ACM-SIAM
Symposium on Discrete Algorithms (SODA 2013), SIAM, 2013, 652--663.

\bibitem{Su71} A.~G.~Sukharev, Optimal strategies of the search for an extremum, USSR Comput.\ Math.\ Math.\ Phys.\
11 (1971), no.~4, 119--137.

\bibitem{TWW} J.~F.~Traub, G.~W.~Wasilkowski, H.~Wo\'zniakowski, Information-Based Complexity, Academic Press,
1988.

\bibitem{Ul26} M.~Ullrich, On bounds between all $s$-numbers and widths of convex sets, preprint, 2026,
arXiv:2608.05024.

\bibitem{Wa84} G.~W.~Wasilkowski, Some nonlinear problems are as easy as the approximation problem, Comput.\ Math.\
Appl.\ 10 (1984), 351--363.

\end{thebibliography}
\end{document}